\documentstyle[amscd,amssymb,verbatim]{amsart}
\newcommand{\lrar}[1]{\begin{picture}(50,10)(-25,-5)
\put(-25,0){\vector(1,0){50}}\put(0,5){\makebox(0,0)[b]{\mbox{$#1$}}}
\end{picture}}
\newcommand{\ldar}[1]{\begin{picture}(10,50)(-5,-25)
\put(0,25){\vector(0,-1){50}}\put(5,0){\mbox{$#1$}}\end{picture}}

\newcommand{\dd}{\partial}

\newcommand{\Tr}{\operatorname{Tr}}
\newcommand{\Nm}{\operatorname{Nm}}

\newcommand{\Gal}{\operatorname{Gal}}

\renewcommand{\Re}{\operatorname{Re}}

\newcommand{\Frob}{\operatorname{Frob}}

\renewcommand{\P}{{\Bbb P}}

\newcommand{\Ind}{\operatorname{Ind}}

\newcommand{\eps}{\epsilon}

\numberwithin{equation}{section}

\newtheorem{thm}{Theorem}[section]
\newtheorem{prop}[thm]{Proposition}
\newtheorem{lem}[thm]{Lemma}
\newtheorem{cor}[thm]{Corollary}
\newenvironment{rem}{\vspace{3mm}\noindent
{\bf Remark.}}{\vspace{3mm}}
\newenvironment{definition}{\vspace{3mm}\noindent
{\bf Definition.}}{\vspace{3mm}}
\newenvironment{rems}{\vspace{3mm}
\noindent {\bf Remarks.}}{\vspace{3mm}}
\newenvironment{ex}{\vspace{3mm}\noindent
{\bf Example.}}{\vspace{3mm}}

\newcommand{\Pf}{\noindent {\it Proof}}

\newcommand{\ov}{\overline}

\newcommand{\ra}{\rightarrow}

\renewcommand{\O}{{\cal O}}

\newcommand{\Res}{\operatorname{Res}}

\newcommand{\la}{\lambda}

\newcommand{\C}{{\Bbb C}}
\newcommand{\R}{{\Bbb R}}
\newcommand{\Z}{{\Bbb Z}}

\newcommand{\Ga}{\Gamma}

\newcommand{\sub}{\subset}
\newcommand{\ed}{\qed\vspace{3mm}}

\title{When is the Fourier transform of an elementary function 
elementary?}
\author{Pavel Etingof, David Kazhdan, and Alexander Polishchuk}

\begin{document}
\maketitle

\bigskip

\centerline{ABSTRACT}
\vskip .1in

Let $F$ be a local field, $\psi$ a 
nontrivial unitary additive 
character of $F$, and $V$ a finite dimensional 
vector space over $F$.
Let us say that a complex function on $V$ is elementary if
it has the form
$$
g(x)=C\psi(Q(x))\prod_{j=1}^k\chi_j(P_j(x)), x\in V, 
$$
where $C\in \C$, $Q$ is a rational function (the phase function), 
$P_j$ are polynomials, and 
$\chi_j$ multiplicative characters of $F$.
For generic $\chi_j$, this 
function canonically extends to a distribution on $V$
(if char($F$)=0).  

Occasionally, the Fourier transform of an elementary function  
is also an elementary function
(the basic example is the Gaussian integral: $k=0$, $Q$ is a 
nondegenerate quadratic form). 
It is interesting to determine when exactly this happens.
This question is the main subject of our study. 

In the first part of this paper we show that for $F=\R$ or $\C$, if 
the Fourier transform of an elementary 
function $g\ne 0$ with phase function $-Q$ such 
that $\text{det}\ d^2Q\ne 0$ is another elementary 
function $g^*$ with phase function $Q^*$, 
then $Q^*$ is the Legendre transform of $Q$ 
(the ``semiclassical condition''). 
We study properties and examples of phase functions satisfying this condition,
and give a classification of phase functions such that
both $Q$ and $Q^*$ are of the form $f(x)/t$, 
where $f$ is a homogeneous cubic polynomial and $t$ 
is an additional variable (this is one of the simplest 
possible situations). Unexpectedly, the proof  
uses Zak's classification theorem for Severi varieties.
\footnote{Unfortunately, this proof turned out to be incomplete.
A complete (different) proof is given in \cite{CS}, Corollary 4.}   
 
In the second part of the paper we give a necessary and sufficient condition 
for an elementary function to have an elementary Fourier 
transform (in an appropriate ``weak'' sense)
and explicit formulas for such Fourier 
transforms in the case when $Q$ and $P_j$ are 
monomials, over any local field $F$. We also describe a generalization of 
these results to the case of monomials of norms 
of finite extensions of $F$. Finally, we generalize 
some of the above results (including Fourier integration
formulas) to the case  
when $F=\C$ and $Q$ comes from a prehomogeneous vector space.

\section{Introduction}

\subsection{Motivations}

Let $F$ be a local field, $\psi :F\rightarrow \C^{*}$ a nontrivial
unitary additive character,  
$V$ a finite dimensional vector space over $F$, and 
$Q$ a nondegenerate quadratic form on $V$.
 It is well known that 
the Fourier transform of the function $\psi (-Q)$ has the form
\begin{equation}
\widehat{\psi (-Q)}=\epsilon (Q)\psi (Q^{-1})\label{Gauss}
\end{equation}
where
$Q^{-1}$ is 
the inverse quadratic form on the dual space $V^*$
(i.e., $dQ\circ dQ^{-1}=Id$). As was shown in
\cite{K1}, Proposition 3, 
there exists an analog of (\ref{Gauss}) for some homogeneous rational
functions on $V$ of homogeneity degree 2 
which are not quadratic polynomials. More
precisely, let  $E$ be a cyclic cubic extension of $F$, and
${\mathcal E}:F^*\rightarrow
\C^{*}$ a nontrivial cubic character which is trivial on the image of the
norm map $\Nm:E^{*}\rightarrow F^{*}$. Let $\phi _{\mathcal E}$ be the
distribution on the vector space $F\oplus E$ such that 
$$\phi _{\mathcal E}(t,x)={\mathcal E} (t)|t|^{-1}\psi
(\Nm(x)/t)$$ Then  we have 
\begin{equation}
\widehat{\phi _{\mathcal E}}=\epsilon \phi _{\mathcal E},\label{cub}
\end{equation}
where $\epsilon$ is $\pm 1$. The proof of
(\ref{cub}) given in \cite{K1} 
is based on the  analysis of the smallest  special
representation 
of the group $D_4 (F)$ and uses global arguments (such as  the
existence of Eisenstein series). We were interested to see whether 
(\ref{cub}) could be proved by local methods and whether there
exist interesting generalizations of 
(\ref{cub}). More precisely, let us say
that a distribution $g$ on a vector space $V$ is "elementary" if
it has the form $C\psi(Q)\chi _1 (P_1)....\chi _k(P_k)$
for some rational function $Q$ (called the phase function), 
polynomials $P_1 ,...,P_k$ on $V$, 
and multiplicative characters $\chi _1,...,\chi _k$ on $F$.
We say that a distribution  $g$ is "special" if both
$g$ and its Fourier transform  $ \hat g$ are "elementary".
We were interested in finding which distributions are ``special''. 

In this paper we present a number of  "special"
distributions. Moreover, we 
show that the description of such distributions is almost
independent of the local field $F$, and the
Fourier transform $ \hat g$, up to a factor, is described
algebraically in terms of  $g$. Therefore this paper provides 
supporting evidence 
for the conjecture of the existence of an 
algebro-geometric integration theory proposed in \cite{K2}.

\subsection{ Statement of the problem.}\label{stat} 
Let us formulate our main question more precisely. 
Keeping the notation of the previous section, set
$$
G_{\chi_1,...,\chi_k}^{P_1,...,P_k,Q}(x,\psi)=
\psi( Q(x))\prod_{j=1}^k\chi_j(P_j(x)).
$$ 
The function $G_{\chi_1,...,\chi_k}^{P_1,...,P_k,Q}$ 
is not always defined on the whole space $V$. 
However, if $F$ has characteristic zero
(and conjecturedly, also in characteristic $p$), 
this function defines a distribution 
on $V$ which meromorphically depends on the characters $\chi_j$
(this follows from the resolution of singularities, or from the
theory of D-modules in the archimedean case). 
In particular, for generic values of $\chi_j$ this function
canonically extends to a distribution on $V$.  

{\bf Question:} For which 
$Q,Q^*,P_i,P_j^*,\chi_i,\chi_j^*$ does one have 
\begin{equation}
\widehat{G_{\chi_1,...,\chi_k}^{P_1,...,P_k,-Q}}=C_\psi
G_{\chi^*_1,...,\chi^*_l}^{P^*_1,...,P^*_l,Q^*},\forall \psi,\label{eq1}
\end{equation}
as distributions on $V^*$?

\subsection{Results of the paper.}

The main results of the paper are as follows. 

In Section 2, using the formal stationary phase method, 
we give a necessary condition for (\ref{eq1}), 
when $Q$ has nonzero Hessian. 
This condition
says that $Q^*$ is the Legendre transform of $Q$. 
We refer to this condition as
``the semiclassical condition'', commemorating the fact 
that identity (\ref{eq1}) can be regarded as a ``quantum mechanical''
formula, from which this condition is deduced by using 
the ``semiclassical'' (i.e., stationary phase) approximation. 
 
In Section 3, we discuss properties and examples 
of phase functions satisfying the semiclassical condition.
The most interesting examples we know come from prehomogeneous vector spaces. 
We classify 
phase functions such that 
both $Q$ and $Q^*$ are of the form $f(x)/t$, where $f$ is a homogeneous 
cubic polynomial and $t$ an additional variable (Section 3). 
The classification says that in this case $f$ is a relative invariant of a  
regular prehomogeneous 
vector space of degree 3 (there are seven cases). 
The proof of this classification theorem is based on 
 Zak's classification theorem for Severi varieties. 

In Section 4, we consider 
elementary functions in which the polynomials 
$P_j$ and the phase function $Q$ are monomial, and
find a necessary and sufficient condition for the Fourier transform of 
such a function to be a function of the same type
(in the ``weak'' sense). 
This condition is an identity with $\Gamma$-functions. 
We generalize this result to the case when 
$Q,P_j$ are monomials of norms of finite extensions of $F$, 
and, more generally, when they are monomials of relative
invariants of prehomogeneous vector spaces over $F$.  

In Section 5, we show that in 
the archimedean case, the condition of Section 4 
for the existence of an integral identity
 can be reformulated in combinatorial terms
(more precisely, in terms of so-called {\it exact covering systems}). 
We write down the explicit integral identities 
in the case when these combinatorial conditions are satisfied. 
We give the simplest nontrivial examples of integral identities,
including the case of prehomogeneous vector spaces.  

In Section 6, we generalize the results of Section 5 
to non-archimedean fields $F$.
Using the known formula for the Gamma function of a cyclic 
field extension, we obtain integral formulas of type (\ref{eq1}).
In particular, we give a new proof of formula (\ref{cub}). 
 
These results have natural analogues in the case when $F$ is a finite field
which will be described in a separate paper by D.K. and A.P. 
The role of distributions in this case is played by perverse sheaves and
the Fourier transform is replaced by its geometric analogue defined by Deligne.
Applying the trace of the Frobenius to identities with perverse
sheaves, one obtains nontrivial elementary identities with
exponential sums.   

\subsection{ Acknowledgements.} We are grateful to D. Arinkin, 
and M. Kontsevich for useful discussions, and to B. Gross for calling our 
attention to Zak's theorem. The work of P.E. and D.K. was
supported by the NSF grant DMS-9700477. The work of A.P. 
was supported by the NSF grant DMS-9700458.

\section{The semiclassical condition}

\subsection{ Formulation of the semiclassical condition.}

Recall the definition of the Legendre transform
(see e.g., \cite{Ar}).
Let $V$ be a finite dimensional real vector space, $v_0\in V$, and 
 $Q$  a smooth function on a neighborhood of $v_0$
such that $\text{det}Q''(v_0)\ne 0$. 
Let $Q'(v_0)=p_0\in V^*$ (where $Q',Q''$ are the first and second 
differentials of $Q$). 
Then the Legendre transform of $Q$ is the 
smooth function $L(Q)$ defined in a neighborhood 
of $p_0$ by $L(Q)(p)=pv_p-Q(v_p)$, where $v_p$ is 
the unique critical point of $pv-Q(v)$
in a neighborhood of $v_0$. 

This definition generalizes tautologically to 
the case when $V$ is a vector space over any field, and 
$Q$ is a regular function on the formal neighborhood of $v_0\in V$. 

It is obvious that if $Q$ 
is an algebraic function then so is $L(Q)$. 

Recall from Section \ref{stat} the definition of 
the function $G_{\chi_1,...,\chi_k}^{P_1,...,P_k,Q}(x,\psi)$. 
For convenience we will always assume that 
the pole divisor of
$Q$ is contained in the divisor $P_1...P_k=0$ (this does not
cause a loss of generality).  

Using the stationary phase method, we will 
prove the following theorem: 

\begin{thm}\label{classical}   
Let $F$ be $\R$ or $\C$. 
Suppose that (\ref{eq1}) is satisfied, and $Q$ has nonzero Hessian. 
Then 

(i) The rational map of algebraic varieties 
$\underline{V}\to \underline{V}^*$
given by $x\to Q'(x)$ is a birational isomorphism. 

(ii) $Q^*$ is the Legendre transform of $Q$. 
\end{thm}

We will call this necessary condition of (\ref{eq1}) {\bf the semiclassical
condition}.

\begin{rem} We expect that a similar result holds 
over non-archimedean local fields. 
\end{rem}

Theorem \ref{classical} is proved in Section 2.3. In the next section, 
we explain the formal stationary phase method, which is necessary for 
the proof. 

\subsection{  The formal stationary phase method.}

The idea of the classical 
stationary phase method can be
summarized as follows. 

Let $V$ be a real finite dimensional vector space with 
a volume form. 
Let $\phi$ be a function defined in an open 
set $B$ around $0$ in $V$ which 
has a nondegenerate critical point at $0$. Let $f$ be a
smooth real-valued function whose support is a compact 
subset of $B$. Consider the integral 
$$
I(\hbar)=\int f(x)e^{-i\phi(x)/\hbar}dx, \hbar>0.
$$

\begin{thm}\label{stphase} (see \cite{AVG} and references therein) 
The function $I(\hbar)$ has the following asymptotic expansion
as $\hbar\to 0$:
$$
I(\hbar)\sim C\hbar^{dim(V)/2}|\text{det}(\phi''(0))|^{-1/2}
e^{-i\phi(0)/\hbar}(f(0)+\sum_{j=1}^\infty R_j(f,\phi)(i\hbar)^j),
$$
where $R_j(f,\phi)=\frac{\hat R_j(f,\phi)(0)}{\text{det}(\phi''(0))^{N_j}}$,
$N_j\in {\Bbb Z}_+$, and $\hat R_j$ are differential polynomials 
with rational coefficients. 
\end{thm}

\begin{rem} The functions $R_j$ are complicated, but there is an algorithm 
of computing them which can be expressed in terms of Feynman diagrams.    
\end{rem}

Now let $p$ be a variable taking values in $V^*$. 
If $p$ is small enough, the function $\phi(v)-pv$ has a unique critical 
point $v_p$ near zero, which is nondegenerate.
Therefore, we have 
$$
\int f(x)e^{i(px-\phi(x))/\hbar}dx\sim 
$$
$$
C\hbar^{dim(V)/2}|\text{det}(\phi_p''(0))|^{-1/2}
e^{-i\phi_p(0)/\hbar}(f_p(0)+\sum_{j=1}^\infty R_j(f_p,\phi_p)(i\hbar)^j),
$$
where $\phi_p(v)=\phi(v+v_p)-p(v+v_p), f_p(v)=f(v+v_p)$. 

Note that $\phi_p(0)=-L(\phi)(p)$, where $L$ is the Legendre transform. 

Now we will generalize this to the formal setting.
Let $V$ be a finite dimensional vector space over a field $F$
of characteristic zero.  
For any regular functions $f,\phi$ on a formal neighborhood of zero 
in $V$, such
that $\phi(0)=0,\phi'(0)=0$, $\text{det}(\phi''(0))\ne 0$,
define the regular
function $J_{f,\phi}(h,p)$ on the formal neighborhood of zero 
in $V^*((h))$ by  
$$
J_{f,\phi}(h,p)=\left(\frac{\text{det}(\phi_p''(0))}
{\text{det}(\phi''(0))}\right)^{-1/2}
e^{-L(\phi)(p)/h}(f_p(0)+\sum_{j=1}^\infty R_j(f_p,\phi_p)h^j)
$$
(this is proportional to the right hand side of the stationary phase 
formula, with $h=i\hbar$). 
We no longer
claim that this series gives the asymptotic expansion 
of the integral $\int f(x)e^{(\phi(x)-px)/h}dx$, 
 because this integral is not defined. 
However, we can still claim that the series $J$ satisfies 
the same differential equations as the integral would satisfy 
if it existed. More precisely, we have the following lemma, which will be 
used to prove Theorem \ref{classical}.  

Let $D$ be a differential operator on $V$ with polynomial coefficients 
over $F((h))$, and let $\hat D$ be the operator on $V^*$
obtained from $D$ by the Fourier automorphism $\frac{\dd}{\dd v}\to h^{-1}v$, 
$p\to -h \frac{\dd}{\dd p}$. 

Define the differential polynomial $E_D(f,\phi)$ by 
 $D(fe^{\phi/h})=E_D(f,\phi)e^{\phi/h}$. 

\begin{lem}\label{formstat} One has 
$$
\hat D J_{f,\phi}(h,p)=J_{E_D(f,\phi),\phi}(h,p).
$$
\end{lem}
  
\Pf . The statement is obvious from the stationary phase formula if
$F=\R$ (by integration by parts), 
and $f,\phi$ are expansions of smooth functions such that $f$  
has compact support inside of the domain of $\phi$. 
Since the statement is purely algebraic, 
it holds in general (because an arbitrary jet can be the jet 
 of a function with compact support). 
\ed

\subsection{  Proof of Theorem \ref{classical}.}

We will consider the case $F=\R$; the case of $\C$ is similar. 

Let $Q,P_j,\chi_j$, $j=1,..,k$ be as in Section 1.1. 
Let $Z\subset V_\C$ be the locus of zeros of $P_i$. 
The function 
$g(x):=G_{\chi_1,...,\chi_k}^{P_1,...,P_k,Q}(x,\psi)$ 
is smooth and nonvanishing on $V_\C\setminus Z$ and hence
generates a 1-dimensional local system on 
$V_\C\setminus Z$. Let us extend this local system to 
an irreducible D-module on $V_\C$ and call this extension 
$M_g$. (Note that all our D-modules are algebraic D-modules).

Consider the D-module $M(r)$ generated by the distribution 
$g_r:=(P_1...P_k)^rg$
for sufficiently large $r$. 
Since $M(r)$ is holonomic, it has finite length, 
and so for sufficiently large $r$, the D-module 
$M(r)$ is independent of $r$. Let us denote this D-module by
$M_\infty$. 

It is easy to see that there is an exact sequence of 
D-modules 
$$
0\to K\to M_\infty\to M_g\to 0,
$$
where $K$ is supported on the divisor $P_1...P_k=0$. 
This implies that 
$M_g$ is isomorphic to 
$D_V/I_r$ (for large enough $r$), where
$D_V$ is the algebra of differential operators on $V$, and 
$I_r$ is the left ideal of differential operators which 
annihilate $g_r$ formally (i.e., outside of $Z$).

\begin{prop}\label{genrank}
For generic $\psi$, the rank 
of the Fourier transform of 
$M_g$ is at least the degree $d$ of the map $Q':\underline{V}\to 
\underline{V^*}$. 
\end{prop}

\Pf . Let us change the ground field from 
$\C$ to $K=\overline{\C((\hbar))}$, and set 
formally $\psi(x)=e^{ix/\hbar}$. 
It is enough to prove the claim for this particular $\psi$. 

For this purpose, it is enough to produce, for a generic 
$p_0\in V_\C^*$, a collection of $d$ linearly 
independent solutions of the differential equations $\hat D\phi=0$,
$D\in I_r$,  
in the formal neighborhood of $p_0$.

To produce such solutions, we will use the formal stationary phase 
method. We will pick $p_0$ generically. Then the equation 
$Q'(x)=p_0$ has exactly $d$ distinct solutions $x_1,...,x_d$, and $Q''$ is 
nondegenerate at all these points. 
Define the power series $\phi_i(x)=Q(x+x_i)-Q(x_i)-p_0x,\ f_i(x)=
\prod_j \chi_j(P_j)(P_1...P_k)^r(x+x_i)$. 
Define 
$$
\eta_j(\hbar,p)=J_{f_j,\phi_j}(i\hbar,p).
$$
It follows from Lemma \ref{formstat} that these series 
are indeed solutions of the equations $\hat D(p+p_0)\phi(p)=0$, 
$D\in I_r$. 

It remains to show that the solutions $\eta_i$ are 
linearly independent. 
To do this, 
consider the power series 
$L(\phi_i)(p)$, and look at their second degree terms, which are equal  
to $Q''(x_i)^{-1}(p,p)$. 

\begin{lem}\label{distinct} For generic $p_0$, 
the forms $Q''(x_i)^{-1}$ are distinct. 
\end{lem}

\Pf . The forms 
$Q''(x_i)^{-1}$ are all the values 
at $p=p_0$ of the multivalued algebraic function
$((Q')^{-1})'(p)$, which is the derivative of the function
$(Q')^{-1}(p)$ that has exactly $d$ branches by the definition. 

But we claim that the derivative of any algebraic vector-function
has at least (and hence exactly) as many branches as the function itself. 
Indeed, it is enough to check it for scalar functions $f(z)$ 
of one variable. But in the one variable case, one always has
$f\in \C(z,f')$, since any monic algebraic equation of degree 
$>1$ satisfied by $f$ over 
$\C(z,f')$ can be differentiated to get a monic equation of lower
degree. The lemma is proved. \qed

Now let us make a change of variables $\hbar\to t^2\hbar$, $p\to t p$, 
and let $t$ tend to $0$. Then we have 
$\eta_j\to e^{iQ''(x_j)^{-1}(p,p)/\hbar}$, which are linearly 
independent functions since $Q''(x_j)^{-1}$ are distinct. 
This implies that $\eta_i$ are linearly independent. 
The proposition is proved.
\ed

Let $G_1=G_{\chi_1,...,\chi_k}^{P_1,...,P_k,-Q}$,  
$G_2=G_{\chi^*_1,...,\chi^*_l}^{P^*_1,...,P^*_l,Q^*}$.

\begin{cor}\label{fourier}
Suppose that (\ref{eq1}) holds.
Then for generic $\psi$ the Fourier transform of the D-module $M_{G_1}$ 
is isomorphic to $M_{G_2}$.
\end{cor}

\Pf .
The relation $\hat G_1= C_\psi G_2$ implies that 
the Fourier transform of the D-module generated by 
$G_1$ is the D-module generated by $G_2$. 
Both of these D-modules 
are holonomic and have only one component of the Jordan-Holder
series which has nonzero rank, namely, $M_{G_1}$ and $M_{G_2}$. 
But by Proposition \ref{genrank}, the rank 
of the Fourier transform of $M_{G_1}$ is positive. 
This implies the statement. 
\ed

Now let us prove the theorem. We start with statement (i). 
By Proposition \ref{genrank}, for generic $\psi$
the rank of the Fourier transform $M_{G_1}$ is at least $d$. 
On the other hand, by Corollary 
\ref{fourier} this Fourier transform is $M_{G_2}$, so its 
rank is $1$. So $d=1$, as desired. 

Statement (ii) follows from the fact that $G_2$ satisfies 
(outside of $Z$) the differential 
equations $\hat D f=0$, where $D\in I_r$, and hence (since 
the rank of the Fourier transform of $M_{G_1}$ is $1$),
the formal expansion of $G_2$ must coincide with the formal series 
constructed using the stationary phase method. This implies that 
$Q^*=L(Q)$, since it is clear that if $J_{f,\phi}(h,p)
=e^{-\xi(p)/h}(1+O(p,h))$
then $\xi=L(\phi)$ (this follows from 
theorem \ref{stphase}). The theorem
is proved.

\begin{rem} It is known that Theorem \ref{stphase} 
has a non-archimedean analogue, 
which is even simpler than this theorem itself: in this case, the asymptotic
expansion of the integral 
contains only the leading term and no higher terms. We expect that this
result can be used to prove Theorem 2.1 in the non-archimedean case.  
\end{rem}

\subsection{ Integral identities in the weak sense and the
  semiclassical condition.}

Let $F$ be archimedean, and $V$ a finite dimensional vector space
over $F$. Let $P$ be a polynomial on $V$ and 
$R$ a polynomial on $V^*$.  
Let $N$ a positive integer, and 
${\mathcal S}_N^{P,R}(V)$ be the space of
Schwartz functions on $V$ of the form $|P|^{2N}f$
 (where $f$ is a
Schwartz function), whose Fourier transform has the form
$|R|^{2N}g$ (where $g$ is a Schwartz function).  
It is easy to construct examples of elements of 
this space: for instance, one can take the function 
$|R|^{2N}(\dd)|P|^{2N'}(x)f$, where $N'>>N$, and $f$ 
is any Schwartz function. 

As we mentioned before, the function
$G^{P_1,...,P_k,Q}_{\chi_1,...,\chi_k}$
defines a distribution only for generic values of the characters
$\chi_j$. However, for any characters $\chi_j$ this function
defines a linear functional on the space 
${\mathcal S}_N^{P_1...P_k,1}(V)$
for large enough $N$. 

\begin{definition} We will say that the integral identity (1.3)
holds {\bf in the weak sense} if it holds on the space 
${\mathcal S}_N^{P_1^*...P_l^*,P_1...P_k}(V^*)$ for large enough $N$. 
\end{definition}

\begin{thm} Theorem 2.1 remains valid if (1.3) holds 
only in the weak sense. 
\end{thm} 

The proof is analogous to the proof of Theorem 2.1.

\section{Rational functions satisfying the semiclassical condition.}

In this section we would like to study systematically the question: 
which rational functions satisfy the semiclassical condition?

\subsection{  Properties of functions satisfying 
the semiclassical condition.}
 
Let $V$ be a finite dimensional vector space over an 
algebraically closed field $F$. Denote by $SC(V)$ the set of 
rational functions $Q$ on $V^*$ 
satisfying the semiclassical condition, i.e., such that 
$Q':V\to V^*$ is a birational isomorphism. 

One can characterize 
elements of $SC(V)$ using the notion of Legendre transform, as follows.

\begin{prop}\label{legrat} A function $Q\in F(V)$
belongs to $SC(V)$
if and only if the Legendre transform of $Q$ is rational.
\end{prop}

\Pf . The ``only if'' part is obvious. To prove the ``if'' part,
let $L(Q)=G$ 
and let us differentiate the equation 
$xQ'(x)-Q(x)=G(Q'(x))$. We get $xQ''(x)=G'(Q'(x))Q''(x)$. 
Since $Q''$ is generically nondegenerate, we get $x=G'(Q'(x))$.
Thus, $G'$ is the inverse to $Q'$. \ed

Nondegenerate quadratic forms  
are the simplest 
examples of elements of $SC(V)$. In the following sections we will 
construct other examples
of elements of $SC(V)$.

\subsection{  The projective semiclassical condition.} $\ $

\begin{definition} A homogeneous rational function $f$ on $V$
is said to satisfy {\it 
the projective semiclassical condition} if the map $x\to f'(x)$ defines 
a birational isomorphism $\P V\to \P V^*$. 
\end{definition}

Denote the set of functions satisfying the projective 
semiclassical condition by $PSC(V)$. 

The relationship between $SC(V)$ and $PSC(V)$, which motivates 
the introduction of $PSC(V)$, is given by the following 
easy lemma. 

\begin{lem}\label{d02} 
Let $g:V\to W$ be a homogeneous of degree $d$ rational map of 
finite dimensional vector spaces, 
which defines a birational isomorphism $\bar g: \P V\to \P W$. 
Then $g$ itself is a birational isomorphism if and only if $d=\pm 1$.
In particular, an element $Q\in PSC(V)$ belongs to $SC(V)$ if and 
only if its homogeneity degree is $0$ or $2$.  
\end{lem}

\Pf . It is enough to prove the first statement; the second statement 
is a special case of the first one for $g=Q'$. 
 
{\bf If.} The condition of the lemma implies that 
for a generic vector $v$ one has  $v=tR(g(v))$, where $R$ is a rational 
function, and $t$ is a factor to be determined. Thus we have 
$t^{d}g(R(g(v)))=g(v)$, which allows one to determine $t$ rationally 
since $d=\pm 1$. 

{\bf Only if.} This part is clear, since a homogeneous birational isomorphism 
between vector spaces
has to have homogeneity degree $1$ or $-1$.
\ed
 
\begin{cor}\label{fproverf} A homogeneous function $f\in F(V)$  belongs to 
$PSC(V)$ if and only if the 
map $f'/f:V\to V^*$ is a birational isomorphism. 
\end{cor}

\Pf . The ``if'' part is clear. The ``only if'' part 
follows when one applies Lemma \ref{d02} to
$g=f'/f$. 
\ed

\begin{cor}\label{constr}
(i) Any nonzero (in $F$) integer power of 
a function $f\in PSC(V)$ belongs to $PSC(V)$. 

(ii) Let $V,W$ be finite dimensional $F$-vector
spaces, and $f\in PSC(V)$, 
$g\in PSC(W)$. 
Then the exterior tensor product $(f\otimes g)(v,w)=f(v)g(w)$ on $V\oplus W$ 
belongs to $PSC(V\oplus W)$.  

(iii)
Any function of the form 
$f_1^{n_1}...f_k^{n_k}$, where $f_i\in PSC(V_i)$, 
and $n_i\in\Z$ are nonzero in $F$, belongs to
$PSC(\oplus_i V_i)$; it belongs to $SC(\oplus_i V_i)$ iff 
$\sum n_id_i=0$ or $2$, where $d_i$ are the homogeneity degrees of $f_i$. 
\end{cor}

\Pf .
Statements (i) and (ii) are obvious from Corollary \ref{fproverf}. 
Statement (iii) follows from (i),(ii), and Lemma \ref{d02}.
\ed

Part (iii) of Corollary \ref{constr} 
allows one to obtain functions 
satisfying the semiclassical condition from functions
satisfying the projective semiclassical condition. 
The simplest example of functions so obtained 
are monomial functions
$x_1^{n_1}...x_k^{n_k}$, where $\sum n_i=0$ or $2$. 

So we will now study the projective semiclassical consition 
more systematically. 

\subsection{ The multiplicative Legendre transform.}

To think about elements of $PSC(V)$, it is useful to introduce
the notion of the multiplicative Legendre transform. 

Let $f$ be a homogeneous function, 
and $det((f'/f)')$ is not identically zero. 
In this case we can define a function 
$f_*$ by $f_*(f'/f(x))=1/f(x)$
(as the usual Legendre transform, it can be defined 
in an analytic as well as a formal setting). 
If $f$ is homogeneous of degree 
$d$ then so is $f_*$. 

\begin{definition} We will call $f_*$ the multiplicative Legendre 
transform of $f$. 
\end{definition}

\begin{rem} Our terminology is motivated by the fact 
that $f_*$ is the multiplicative Legendre 
transform of $f$ if and only if $L(\ln f)=d+\ln f_*$
(over $\C$). 
\end{rem}

It is obvious that the operation $f\to f_*$ commutes with 
exterior tensor product. Also, $(f^n)_*=n^{-nd}f_*^n$, 
 where $d$ is the degree of $f$ (if $n\ne 0$ in $F$).

\begin{ex} $(\prod x_i^{n_i})_*=\prod n_i^{-n_i}\prod
x_i^{n_i}$. 
\end{ex}

\begin{prop}\label{mlt} 
(i) $f_*'/f_*\circ f'/f=Id$. 

(ii) $f_{**}=f$. 
\end{prop}

\Pf . The first identity is obtained by differentiating the definition 
of $f_*$. The second one is obtained by applying $f_{**}$ to both sides of 
(i). 
\ed

One can characterize 
elements of $PSC(V)$ using the notion of 
the multiplicative Legendre transform, as follows.

\begin{prop}\label{mltrat}  A homogeneous 
rational function $f$ belongs to $PSC(V)$ 
if and only if $f_*$ is rational. 
\end{prop}

\Pf . The ``only if'' part follows from Corollary \ref{fproverf}. 
The ``if'' part follows from Proposition \ref{mlt} (i).  
\ed

\begin{rem} We see from the above that
elements of $PSC(V)$ are multiplicative analogs of 
elements of $SC(V)$. 
\end{rem}

\begin{ex} Let us point out an easy method of creating new 
functions satisfying the projective semiclassical condition
out of ones already known. It is straightforward to compute that if
$f_*$ is the multiplicative Legendre transform of $f$ on $V$ then 
the multiplicative Legendre transform of the function 
${\Bbb F}(x,y)=f'(x)y+f(x)$ on 
$V^2$ is ${\Bbb F}_*(x_*,y_*)=
(d-1)^{1-d}(f_*'(y_*)x_*-f_*(x_*))^{d-1}f_*(y_*)^{2-d}$. 
This formula is valid also for $d=1$ if we agree that $0^0=1$. 
\end{ex}

\begin{rem} Note that if ${\Bbb F}$ is a polynomial 
and $d\ge 3$ then ${\Bbb F}_*$ is not a polynomial. 
Construction of polynomial elements $f\in PSC(V)$ such that 
$f_*$ is also a polynomial is more tricky, and the only examples we know  
are described in the next section. 
\end{rem}

\subsection{  Construction of elements of $PSC(V)$
from prehomogeneous vector spaces.}

Recall \cite{Sa,KS} that a prehomogeneous vector space over $F$ is a triple 
$(G,V,\chi)$, where $G$ is an algebraic group over $F$, 
$V$ an algebraic representation of $G$, and $\chi$ a nontrivial 
algebraic character 
of $G$ such that 

(i) $V$ has a Zariski dense $G$-orbit, 
and 

(ii) there exists a nonzero polynomial $f$ on $V$ such that 
$f(gv)=\chi(g)f(v)$, $g\in G$ (it is obvious that if condition (i) holds, 
$f$ is unique up to a scalar). 

Here we will assume that the group $G$ is reductive. 

\begin{rem} Prehomogeneous vector spaces 
were introduced in the 60-s by Sato (see \cite{Sa}). 
Prehomogeneous vector spaces over $\C$ with reductive $G$ 
and irreducible representation $V$ 
have been classified, see \cite{KS}.
\end{rem}

Let $G_0=Ker(\chi)$. Then $G_0$ is a codimension $1$ subgroup of $G$, 
and the function $f$ is invariant under $G_0$.
The following proposition characterizes the ring of all 
$G_0$-invariants. 

\begin{prop}\label{ringinv} 
\cite{KS} The ring of $G_0$-invariants of $V$ is 
$F[f]$. 
\end{prop}

\begin{cor}\label{dense} $G_0$ has a dense orbit 
in $\P V$. 
\end{cor}

\Pf . (of Proposition \ref{ringinv})
 
Let $h$ be any homogeneous polynomial which is an eigenfunction 
for $G$ and invariant under $G_0$. Since $G/G_0={\Bbb G}_m$, 
there exist nonzero integers $k,l$ such that $h^kf^l$ is $G$-invariant. 
Since the $G$-action has a dense orbit, this implies that 
$h^kf^l$ is a constant, i.e., $h,f$ 
are powers of the same polynomial $g$, which is invariant under $G_0$ 
and is an eigenfunction of $G$. By the definition of $f$, this polynomial 
has to be $f$ (we can't have $g^k=f$ for $k>1$ since then there would exist
elements of $G$ that are not in $G_0$ but preserve $f$). 
The proposition is proved.  
\ed

The function $f$ (which is uniquely determined up to scaling) 
is called the relative invariant of $(G,V,\chi)$. 

 From now till the end of the subsection
we will assume that the characteristic of $F$ is zero.

Recall \cite{Sa,KS} that a regular prehomogeneous vector space 
is such that $det(f'')$ is not identically zero. 

\begin{prop}\label{regpreh} (see \cite{Sa}) Let $(G,V,\chi)$ be a 
regular prehomogeneous vector 
space, and $f$ its relative invariant. Then $f\in PSC(V)$,
and its multiplicative Legendre transform is a polynomial.
\end{prop}

\Pf . It is clear that $(G,V^*,\chi^{-1})$ is a prehomogeneous vector 
space. Indeed, the existence of an open orbit, and 
the existence of a relative invariant of the same degree $d$ 
follows from the fact that the representation $S^kV^*$ is completely 
reducible. 

Let $f_*$ be the relative invariant of $(G,V^*,\chi^{-1})$. 
Consider the function $f_*\circ f'$ on $V$. This function 
is nonzero because of regularity, has degree $d(d-1)$, and 
is $G_0$-invariant. Thus, this function is proportional 
to $f^{d-1}$. So $f_*$ is the multiplicative Legendre transform of $f$ 
up to scaling. The proposition is proved.    
\qed

\noindent{\bf Examples.} 1. $G=GL(1)$, $V=F$, $\chi(a)=a^n$, $f(x)=x^n$. 

2. $G=GL(1)^{n}$, $V=F^n$, the action is 
$(a_1,...,a_n)(x_1,...,x_n)=(a_1x_1,...,a_nx_n)$, 
$\chi(a)=a_1...a_n$, $f=x_1...x_n$. 
 
3. $G=GL(n)$, $V=S^2F^n$, where $F^n$ is the vector 
representation, $\chi=det^2$, $f=det$.  

4. $G=GL(n)\times GL(n)$, $V=F^n\otimes (F^n)^*$,  
$\chi=det\otimes det^{-1}$, $f=det$. 

5. $G=GL(2n)$, $V=\Lambda^2F^{2n}$, $\chi=det$, 
$f=Pf$ (the Pfaffian).

6. $G=E_6\times GL(1)$, $V$ is the 27-dimensional irreducible representation
(with $GL(1)$ acting by scalar multiplication), 
$\chi$ is $z\to z^3$, 
$f$ is the invariant cubic form. 

7. $G=O(N)\times GL(1)\times GL(1)$, $V=F^N\oplus F$, 
where $F^N$ is the vector representation of $O(N)$, 
the action $(g,a,b)(v,x)=(agv,bx)$, $\chi(g,a,b)=a^2b$, 
$f(v,x)=Q(v)x$, where $Q$ is the invariant quadratic form. 

To conclude this subsection, we would like to raise two questions. 

{\bf Question 1.} Is it true that any polynomial $f\in PSC(V)$
such that $f_*$ is also a polynomial, is a relative invariant
of a prehomogeneous vector space?

{\bf Question 2.} Is it true that any polynomial $f\in PSC(V)$
such that $f_*$ is also a polynomial, is rigid? In other words, 
is the set of equivalence classes of such polynomials finite for each degree? 
 
\begin{rem} 
Clearly, a positive answer to question 1 implies a positive answer 
to question 2, but question 2 could be more tractable. 
\end{rem}

For degree $\le 3$, the answer to both questions is yes. 
This is proved in the next subsection. 
 
\subsection{ Classification of cubic forms with cubic multiplicative Legendre 
transform.}

In this section the field $F$ has characteristic zero. 
We prove the following theorem. 

\begin{thm}\label{cubic} Let $f$ be a cubic form on a finite dimensional 
vector space $V$ such that its multiplicative Legendre transform is also 
a cubic form. Then $f$ is given by one of Examples 1-7 of the 
previous section. 
\end{thm}

The rest of the section is the proof of this theorem. 

({\bf Warning:} It was pointed out by P. Sabatino 
and F. Viviani that this proof is incomplete. 
Namely, the argument with the Hessian 
at the end of the proof of Proposition \ref{316}
is not, by itself, sufficient to conclude that $Z\setminus 0$ is smooth.
However, a different proof of Theorem \ref{cubic} has been given 
in \cite{CS}.) 

First of all, we may assume that $f$ (and $f_*$) are irreducible
(in which case $dim(V)\ge 3$). 
If any of them is reducible, it is a product of a linear and a quadratic 
form or of three linear forms, and it is easy to see that it is given 
by examples 1,2, or 7. 

The functions $f,f_*$ satisfy the equations 
\begin{equation}
f_*(f')=f^2, f(f_*')=f_*^2, f'\circ f_*'(x)=f_*(x)x, 
f_*'\circ f'(x)=f(x)x.\label{defeq}
\end{equation}

Let $X\subset V$ be the zero locus of $f$. 
Let $Z\subset X$ be the zero locus of $(f,f')$, i.e., the singular locus 
of $X$. Define $Z_*\subset X_*\subset V^*$ in a similar way.   
  
\begin{lem}\label{311} (i) $Z=(f')^{-1}(0)$, $Z_*=(f_*')^{-1}(0)$.

(ii) $X=(f')^{-1}(Z_*)$, $X_*=(f_*')^{-1}(Z)$.
\end{lem}

\Pf . Part (i) follows from the fact that $f'=0$ implies $f=0$, 
and similarly for $f_*$ (Euler equation). Part (ii) follows from equations 
(\ref{defeq}). 
\ed

\begin{definition} A 3-dimensional subspace in $V$ is said 
to be a Cremona subspace if it contains three noncoplanar lines
(through $0$) in $Z$ 
and is not entirely contained in $X$. 
\end{definition}

Cremona subspaces in $V^*$ are defined similarly.

\begin{lem}\label{312} Let $L$ be a Cremona 
subspace in $V$. Then 
$f'(L)$ is a Cremona subspace in $V^*$, and the map $f': L\to f'(L)$ 
is the Cremona map $(x,y,z)\to (yz,zx,xy)$ in some coordinates. 
In particular, the intersection of a Cremona subspace with $Z$ 
is the union of the three lines from the definition.  
\end{lem}

\Pf .
Let $z_1,z_2,z_3$ be noncoplanar lines in $Z$ that are contained in $L$, 
and let $p_3=z_1z_2$, $p_2=z_1z_3$, $p_1=z_2z_3$ be the 
planes spanned by pairs of lines. 
Consider the restriction of the map $f'$ to $p_1$. This map is 
given by a homogeneous quadratic form of two variables
with values in $V^*$, 
and it vanishes at lines $z_2,z_3$ by Lemma \ref{311}. 
Therefore, $f'|_{p_1}$ has the form $q(v)w$, where 
$q$ is a quadratic function, and $w\in V$. The same applies to 
$p_2,p_3$. Let $z_1^*,z_2^*,z_3^*$ be the images of
 $p_1,p_2,p_3$. Then $z_i^*\subset Z_*$ ($p_i\subset X$
because $X$ is a cubic) and each $z_i^*$ is a line or zero. 
Let $L_*$ be a 3-subspace in $V^*$ 
containing all $z_i^*$. Let $l$ be a generic plane in $L$. 
Then $f'(L)$ has three intersection  lines with $L_*$ (images of intersections 
of $l$ with $p_1,p_2,p_3$). Since $f'$ is a quadratic map, we have 
$f'(l)\subset L_*$ and thus $f'(L)\subset L_*$. 

Let us show that $z_i^*$ are in fact nonzero (i.e. lines) and noncoplanar. 
If $z_i^*$ 
lie in a 2-plane $p$, $f_*'\circ f'|_L$ is a map from a 3-space to a 2-plane. 
Therefore, $L$ must be contained in $X$, as it is obvious 
from (\ref{defeq}) that $f_*'\circ f'$ is finite (4 to 1) outside of $X$. 
This contradicts the definition of $L$. 

This shows that $L_*$ is a Cremona subspace (it is clear that $L_*$ is 
not entirely contained in $X_*$), and it is obvious that 
if $z_i,z_i^*$ are used as coordinate axes then $f'$ becomes the Cremona map. 
The lemma is proved. 
\ed

\begin{prop}\label{313} A generic plane $p$ (through $0$) in $V$ 
is contained in a Cremona subspace. 
\end{prop}

\Pf . Let $a_1,a_2,a_3$ be the three lines of intersection 
of $p$ with $X$. Let $z_i^*=f'(a_i)\in Z_*$. It is clear that $z_i^*$ 
are noncoplanar, since otherwise $p\subset X$, and we assumed that 
$p$ is generic. Thus the 3-space $L_*$ spanned by $z_i^*$ is a 
Cremona subspace. Let $L=f_*'(L_*)$. Then 
by Lemma \ref{312} $L$ is a Cremona subspace 
which contains $a_i$.
\ed

\begin{prop}\label{314} $dim(Z)\ge \frac{2}{3}dim(V)-1$, 
and $Z$ is not contained in any hyperplane. 
\end{prop}

\Pf . Define a map 
$g:Z\times Z\times Z\times F^3\to V\oplus V$ by 
$G(\zeta_1,\zeta_2,\zeta_3,b_1,b_2,b_3)=
(\sum \zeta_i,\sum b_i\zeta_i)$. This map is 
dominant, since by Proposition \ref{313}
two generic vectors are contained in a Cremona subspace.
This implies the proposition.
\ed

\begin{prop}\label{315}
$Z$ is irreducible. 
\end{prop}

\Pf . It is enough to show that the map $f_*':X_*\to Z$ is surjective, 
i.e., that for any $\zeta\ne 0$, $\zeta\in Z$ there exists 
$x\in X_*$ such that $f_*'(x)=\zeta$. To do this, it is enough to show 
that any such $\zeta$ is contained in a Cremona subspace, since then 
$\zeta$ lies in the image of a special plane in the dual Cremona subspace
in $X_*$.  

Suppose that $\zeta_1,\zeta_2\in Z$ 
are such that $f''(\zeta)(\zeta-\zeta_1,\zeta-\zeta_2)\ne 0$. 
Then \linebreak
$f(\zeta+t(\zeta-\zeta_1)+s(\zeta-\zeta_2))$ is of order 
$ts$ modulo cubic and higher terms 
at $t,s=0$, and thus the points $0,\zeta,\zeta_1,\zeta_2$ 
are not in the same plane and span a Cremona subspace. 

Thus, if the proposition was false, we would have
 $f''(\zeta)(\zeta_1-\zeta,\zeta_2-\zeta)=0$
for all choices of $\zeta_1,\zeta_2$. 
But since $Z$ does not lie in any hyperplane, 
the possible vectors $\zeta-\zeta_1$ span $V$, and similarly for 
$\zeta-\zeta_2$. Thus, $f(\zeta)=f'(\zeta)=f''(\zeta)=0$.
This implies that $f(y+\zeta)=f(y)$ for all $y\in V$, 
and hence $f$ is pulled back from $V/<\zeta>$. This is a contradiction, 
since we assumed that $det(f'')$ is not identically zero. 
\ed    

\begin{prop}\label{316} $dim(Z)=\frac{2}{3}dim(V)-1$
(in particular, $d$ is divisible by $3$), and $Z\setminus 0$ 
is smooth.
\end{prop}

\Pf . Since $Z$ is irreducible, to prove the first statement 
it is sufficient to show that the map $g$ 
defined in the proof of Proposition \ref{314} is generically finite.
For this, it suffices to show that the Cremona subspace containing a generic 
plane is unique. Indeed, if $p$ is a generic plane, $a_1,a_2,a_3$ 
is as in the proof of Proposition \ref{313}, then $a_i$ are contained in any 
Cremona subspace containing $p$, so $f'(a_i)=z_i^*$ are contained  
in the image of this subspace. But $z_i^*$ are noncoplanar, so
such a subspace 
is unique. The first statement is proved. 

Let us now prove that $Z\setminus 0$ is smooth. The dimension of the generic fiber
of the map $f_*':X_*\setminus Z_*\to Z\setminus 0$ is $dim(X_*)-dim(Z)=
d/3$, where $d=dim(V)$. Consider $d f_*': V^*\to V$. It is enough to show that 
the nullity of $df_*'=f_*''$ is $\le d/3$ at all points 
of $X_*\setminus Z_*$ (then it is exactly 
$d/3$ everywhere, and $\pi$ is a smooth fiber bundle). 

Now we will need two simple lemmas. 

{\bf Lemma 1.} $det(f_*'')=const\cdot f^{d/3}$. 

{\it Proof.} This follows from the fact that this determinant 
is nonzero outside of $X$ (where $f_*'$ is a 2-1 covering), the 
irreducibility of $f$, and 
the fact that $\text{deg}(\text{det} f_*'')=d$.

{\bf Lemma 2.}  Suppose that $A(t)$ is a polynomial family of matrices 
such that $det A(t)$ vanishes exactly to the $d$-th order at $t=0$. 
Then the nullity of $A(0)$ is at most $d$.

{\bf Proof.} 
We can assume that the kernel of $A(0)$ 
is the span of the first $r$ basis vectors, so 
the first $r$ columns of $A$ vanish at $0$. 
Thus, $det(A)=O(t^r)$, so $r\le d$.  

Let us now apply Lemma 2 to $f_*''$ restricted to a line transversal to
$X_*$ at a nonsingular point. Taking into account 
Lemma 1, we get that the nullity of $f_*''$ at this point is $\le d/3$,  
which completes the proof of the proposition. 
\footnote{{\bf Warning:} Unfortunately, this 
argument is not sufficient to establish smoothness of $Z\setminus 0$ (we thank P. Sabatino and F. Viviani for pointing this out). 
However, the proposition is valid, as is Theorem \ref{cubic}, as shown in \cite{CS}, Corollary 4
by a different method.}  
\ed

Now we will finish the proof of the theorem. The above proposition 
shows that the projectivization ${\Bbb P}(Z)$ of $Z$ is a 
smooth projective variety of 
dimension $\frac{2}{3}(dim \P (V)-2)$ in $\P (V)$, which does not lie in a 
hyperplane. It has the following property: the line connecting any two points 
of $Z$ is entirely contained in $X$ and thus never goes through a 
generic point of 
$\P (V)$. Varieties with these properties are called Severi
varieties.

Now comes the central part of the proof, which 
is the use of the following classification theorem 
of Severi varieties, due to F. Zak. 
  
\begin{thm} \cite{Z}
Let $Y$ is a smooth, closed subvariety of the 
complex projective space $\C P^{d-1}$, which does not lie in a hyperplane. 
Suppose that 

(i) any line connecting two points of $Y$
belongs to a certain hypersurface, and 

(ii) $\text{dim}(Y)=\frac{2}{3}(d-3)$. 

Then $Y$ is projectively equivalent 
to the singularity locus of the equation 
$f=0$ on $\P (V)$, where $(V,f)$ is one of the following four 
prehomogeneous vector spaces with relative invariant:

1. $V$ is the space of 3 by 3 symmetric matrices, 
$f=\text{det}$. 

2. $V$ is the space of 3 by 3 matrices, $f=\text{det}$.

3. $V$ is the space of skew-symmetric 6 by 6 matrices,
   $f=\text{Pf}$. 

4. $V$ is a 27-dimensional irreducible representation of $E_6$,
$f$ is the invariant cubic form. 
\end{thm}  

This theorem shows that in our situation, 
$(V,f)$ is given by one of the examples 3-6
in the previous section. 

The variety $X$ is reconstructed from $Z$ as the set of points on lines
connecting two points of $Z$. 
The theorem is proved. \qed

\begin{cor}\label{fovert}
Let $Q(x,t)=f(x)/t$, where $f$ is a homogeneous cubic polynomial
on some finite dimensional complex vector space $W$, and $t$ 
is an additional variable. Then the following conditions are equivalent:

(i) $Q\in SC(W\oplus \C)$, and $Q^*=\tilde f(x_*)/t_*$, 
where $\tilde f$ is a cubic polynomial on $W^*$. 

(ii) The polynomial $f$ is as in Examples 1-7.
\end{cor} 

\Pf . It is easy to check directly that $L(f(x)/t)=\tilde f(x_*)/t_*$ 
iff $\tilde f=-f_*$. Thus the statement follows from Theorem \ref{cubic}.
\ed

\begin{rems} 1. In examples 3-5, the map $g$ of Proposition \ref{314} 
has a classical linear algebraic interpretation. Example 3 corresponds 
to simultaneous diagonalization of two quadratic forms. 
Example 4 corresponds to simultaneous diagonalization
of two hermitian forms. Example 5 corresponds to simultaneous
reduction of two skewsymmetric forms to the canonical form
(sum of 2-dimensional forms).  
 
2. Examples 2-7 correspond to the semisimple Jordan algebras of 
degree 3 \cite{J}. It is possible that one can check directly
(i.e., without using Theorem 3.10) that 
in the assumptions of Theorem 3.10, if $dim(V)>1$ then $V$ admits a structure 
of a separable (hence semisimple) Jordan algebra 
of degree 3 such that $f$ is its determinant polynomial.  
This would allow to give another proof of Theorem 3.10 
which would use Albert's theorem
on the classification of simple Jordan algebras, rather than Zak's 
classification theorem.  
\end{rems}

\section{Integral identities with monomials}

\subsection{  Fourier transform.}

In the next four sections we will recall some basic facts 
about analysis over local fields. 
The basic reference for these facts is 
the book \cite{GGP}.

Let $F$ be a local field. We fix a nontrivial additive character $\psi$.  
For any finite dimensional vector space $V$ over $F$, let
${\mathcal S}(V)$ be the space of 
(complex-valued) Schwartz functions on $V$.

For any Haar measure $dx$ on $V$, one  
defines the Fourier transform ${\mathcal S}(V)\to 
{\mathcal S}(V^*)$ by 
\begin{equation}
\hat f(y)=\int_F\psi(yx)f(x)dx.
\end{equation} 

Any Haar measure on $V$ defines a positive inner product on ${\mathcal S}(V)$. 
Let us say that Haar measures $dx$ on $V$ and $dx^*$ on $V^*$ 
are compatible if the Fourier transform is a unitary operator
with respect to this inner product. 
It is easy to see that this condition is symmetric, and that if it is satisfied
then 
one has the inversion formula 
$\hat{\hat f}(x)=f(-x)$. 

If $V$ is identified with $V^*$, one can choose a unique Haar measure 
which is compatible with itself. For example, this is so if 
$V=F$ or $F^n$. From now on, we will use the notation $dx$ 
for this special measure. 

The measure $dx$ on $F$ depends on $\psi$. For example, in the 
archimedean case, if $\psi(x)=e^{i\text{Re}(ax)}$
then $dx$ is $(|a|/2\pi)^{\dim_\R F/2}$
times the Lebesgue measure.  
Below, we use the character $\psi(x)=e^{i\text{Re} (x)}$ 
for $F=\R,\C$, and a character of norm $1$ for the
non-archimedean case. This completely determines $dx$. 

Let ${\mathcal D}(V)$ be the space of distributions on $V$. 
If $V$ carries a Haar measure $dx$, then 
along with the Fourier transform of Schwartz functions,
one can define the Fourier transform of distributions 
${\mathcal D}(V^*)\to {\mathcal D}(V)$ (by duality). 
It will also be denoted by $f\to \hat f$. 
In the case of compatible measures, it satisfies 
the inversion formula.

\subsection{ Multiplicative characters.}

Let $F^*$ be the multiplicative group of $F$, $U(F^*)$  the group 
of continuous unitary characters of $F^*$, and $X(F^*)$  the group 
of continuous 
characters of $F^*$ 
into $\C^*$. If $F=\R$, then $U(F^*)=\R\times \Z/2\Z $ 
and $X(F^*)=\C\times \Z/2\Z $. If $F=\C$, then $U(F^*)=\R\times \Z$, 
and $X(F^*)=\C\times \Z$. In the non-archimedean case,
$U(F^*)$ is $\R/\Z\times D$, where $D$ is a discrete 
countable group, and $X(F^*)=\C/\Z\times D$. 

Let $d_mx$ be the multiplicative Haar measure
on $F^*$, normalized so that $d_mx/dx=1$ at $x=1$. 
Let $\nu_F(x)=(d_mx/dx)^{-1}$ be the norm of $x$
(in the archimedean case, it equals $|x|^{\text{dim}_\R F}$).

For $\la\in X(F^*)$, denote by $\Re \la$ the real number defined by
$|\la(x)|=\nu_F(x)^{\Re \la}$. 
It is obvious that $X(F^*)=U(F^*)\times \R$, 
via $\chi\to (\frac{\chi}{|\chi|},\Re\chi)$. 
We think of the first coordinate as the imaginary part of $\chi$
and of the second as the real part of $\chi$. 

Let us say that $\la\in X(F^*)$ is a singular character
if $\lambda(x)=\nu_F(x)^{-1}$ (in the non-archimedean case), 
$\la(x)=\nu_F(x)^{-1}x^{-n}$, $n\in\Z_{\ge 0}$ (for
$F=\R$), and  $\la(x)=\nu_F(x)^{-1}x^{-n}\bar x^{-m}$,
$n,m\in\Z_{\ge 0}$ (for $F=\C$). It is well known 
that $\la(x)$ is a holomorphic family of distributions on $F$ 
depending on $\la\in X(F^*)$, with simple poles at singular
characters. 

\subsection{ Gamma functions.}

Now define the Gamma function $\Gamma^F(\la)$ of
a local field $F$  to be the meromorphic function on $X(F^*)$ given by 
\begin{equation}
\widehat{\la\nu_F^{-1}}(x)=\Gamma^F(\la)\la^{-1}(x),\label{4.4}
\end{equation}
(whenever both $\lambda\nu_F^{-1}$ and $\lambda^{-1}$ 
define distributions on $F$). 

 From the inversion formula for the Fourier transform 
one gets the functional equation
\begin{equation}
\Gamma^F(\la)\Gamma^F(\nu_F\la^{-1})=\la(-1)\label{4.5}
\end{equation}

Let us give the expressions for the Gamma functions 
of $\R$ and $\C$. 

\begin{lem}\label{gammarc} For $s\in \C$,
let $\la_{s,n}(x)=|x|^s(x/|x|)^n$, $x\in F$. Let 
$\Gamma^F_n(s)=\Gamma^F(\la_{s,n})$. Then 
\begin{equation}
\Gamma^\R_n(s)=(2\pi)^{-1/2}\cdot 2
i^n\Gamma(s)\cos(\pi (s-n)/2),\ n\in \Z/2\Z ,\label{4.6}
\end{equation}
and
\begin{equation}
\Gamma_n^\C(s)=(2\pi)^{-1}
2^{s}i^n\Gamma(\frac{s+n}{2})\Gamma(\frac{s-n}{2})
\sin(\pi (s-n)/2),\ n\in\Z.
\label{4.7}
\end{equation} 
\end{lem}

This lemma is well known and is proved in a straightforward way. 

One can also easily compute the Gamma function 
of a power of $\nu_F$ in the non-archimedean case. 
If $q$ is the order of the residue field of $F$ then 
\begin{equation}
\Gamma^F(\nu_F^s)=\frac{1-q^{s-1}}{1-q^{-s}}.\label{4.7a}
\end{equation}

 From these formulas it is clear that the Gamma function 
has simple poles at the characters $\la\nu_F$, where $\la$ is
singular. It is clear from the definition that it is holomorphic 
everywhere else. (For an exact expression of the Gamma function
in the non-archimedean case, see \cite{GGP}). 

\subsection{  Integral representation of the additive character.}

Let $du$ be the Haar measure on $U(F^*)$ for which 
the Mellin transform $L^2(F^*,d_mx)\to L^2(U(F^*),du)$ is a unitary 
operator. 

We would like to consider a distribution on $U(F^*)$ of the form 
$$
\phi\to \int_{U(F^*)}\Gamma^F(u)\phi(u)du.
$$
Since $\Gamma^F(u)$ has a pole at 
the trivial character, it is necessary to choose 
a regularization of this integral. Our choice, here 
and throughout, will be the following: 
to avoid the pole, the contour of integration 
in the connected component of the identity in 
$X(F^*)$ should be indented in the direction 
of positive values of $\text{Re}(u)$. The following lemma shows that
this choice coincides with the choice of \cite{GGP}, where 
$\Gamma^F(u)$ is defined as the Mellin transform of $\psi$. 
 
\begin{lem}\label{intrep} 
The distribution $\psi(x)$ on $F$ has the following integral
representation:  
$$
\psi(x)=\int_{U(F^*)}\Gamma^F(u)u^{-1}(x)du.
$$
\end{lem}

\begin{rem} This integral is divergent for any 
concrete value of $x$ but is absolutely convergent on 
any test function from the Schwartz space. 
\end{rem}

\Pf . It is easy to see from the definition of the Gamma function that
the Fourier transform on the group $F^*$ of the distribution $\psi(x)$  
(i.e., the Mellin transform)
is equal to $\Gamma^F$ outside of $0$. 
Therefore the statement of the lemma on test functions which vanish at 0
is obtained by applying the inversion 
formula for the Fourier transform. It suffices now to check 
this identity on one test function which does not vanish at 0. 
In the non-archimedean case, this is easy to do for the 
characteristic function of the integers, and 
in the archimedean case one can do it for the function $e^{-|x|^2/2}$.
\qed

\subsection{ The generalized Gamma function.}

For a positive integer $d$,
define a meromorphic function 
$\Gamma^F_{d,a}$ on $X(F^*)$ by the formula
$$
\Gamma^F_{d,a}(\la)=d^{-1}\sum_{\mu:\mu^d=\la}\Gamma^F(\mu)\mu^{-1}(a).
$$

We have $\Gamma_{1,a}^F(\la)=\Gamma^F(\la)\la^{-1}(a)$, so 
the function $\Gamma_{d,a}^F$ is a generalization of the Gamma function. 
We will call $\Gamma^F_{d,a}$ the generalized Gamma function. 
 
  Let us compute the generalized Gamma function in the archimedean case. 

  If $F=\C$ then 
\begin{equation} 
\Gamma_{d,a}^F(\la_{s,nd})=\Gamma^F(\la_{s/d,n})\la_{-s/d,-n}(a),
\end{equation}

  If $F=\R$ and $d=2k+1$ then 
\begin{equation} 
\Gamma_{d,a}^F(\la_{s,n})=\frac{1}{d}\Gamma^F(\la_{s/d,n})\la_{-s/d,-n}(a).
\end{equation}

 If $F=\R$, $\psi(x)=e^{ix}$, $a>0$, and $d=2k$ then 
\begin{equation}
\Gamma_{d,\pm a}^F(\la_{s,0})=(2\pi)^{-1/2}k^{-1}
e^{\pm \pi is/4k}\Gamma(s/2k)a^{-s/2k}.
\end{equation}

\subsection{  The distribution $G_{\la_1,...,\la_k,a}^{n_1,...,n_k}$ on 
$F^k$ and its integral representation.}

 Let $n_1,...,n_k$ be integers, and
$\la_1,...,\la_k\in X(F^*)$. Consider the function on 
$(F^*)^k$ defined by 
\begin{equation}
G^{n_1,..,n_k,a}_{\la_1,...,\la_k}(x_1,...,x_k)=
\psi(a\prod_{i=1}^kx_i^{n_i})\la_1(x_1)...\la_k(x_k).\label{4.10}
\end{equation}
(Here $a\in F^*$ is a parameter). 

$G^{n_1,..,n_k,a}_{\la_1,...,\la_k}$ is a distribution on $F^k$ 
which is holomorphic in $\la_i$ if $\Re \la_i>-1$. 
 
\begin{lem}
For $\la_i$ with real parts $>-1$ one has 
$$
G^{n_1,..,n_k,a}_{\la_1,...,\la_k}(x_1,...,x_k)=
\int_{U(F^*)}\Gamma^F(u)u^{-1}(a)\la_1u^{-n_1}(x_1)...\la_ku^{-n_k}(x_k)du.
$$
(in the sense of distributions).
\end{lem}

\Pf . Follows directly from Lemma \ref{intrep}.
\qed

It will be convenient for us to understand the 
function $G_{\lambda_1,...,\lambda_k}^{n_1,...,n_k,a}$ as a 
distribution {\bf in the
weak sense}. Namely, 
for polynomials $P,R$ define the space ${\mathcal S}_N^{P,R}(V)$
as in section 2 for the archimedean case, and as the space of
Schwartz functions vanishing on the variety $P=0$ 
whose Fourier transform vanishes on $R=0$, 
in the non-archimedean case (so in the non-archimedean case 
it is independent of $N$). As in Section 2, 
for any $\lambda_i$ the function 
$G_{\la_1,..,\la_k}^{n_1,...,n_k,a}$ defines a linear functional on
the space ${\mathcal S}_N^{P_1...P_k}(F^k)$ for a large enough
$N$. We call such a functional a distribution in the weak sense. 

The following lemma is a straghtforward generalization 
of the previous lemma. 

\begin{lem}
For any $\la_i$ one has 
$$
G^{n_1,..,n_k,a}_{\la_1,...,\la_k}(x_1,...,x_k)=
\int_{U(F^*)}\Gamma^F(u)u^{-1}(a)\la_1u^{-n_1}(x_1)...\la_ku^{-n_k}(x_k)du.
$$
(as distributions in the weak sense).
\end{lem}  

\subsection{  The Fourier transform of 
the distribution $G^{n_1,..,n_k,a}_{\la_1,...,\la_k}$.}

\begin{lem} One has 
\begin{equation}
\begin{align*}
&\widehat {G^{n_1,..,n_k,a}_{\la_1,...,\la_k}}(p_1,...,p_k)=\\
&\int_{U(F^*)}\Gamma^F(u)u^{-1}(a)\prod_{i=1}^k\Gamma^F(\la_iu^{-n_i}\nu_F)...
\la_i^{-1}u^{n_i}\nu_F^{-1}(p_i)du.
\end{align*}
\end{equation}
\end{lem}

\Pf . This follows from the previous lemma.
\qed

\subsection{  Identities with monomials.}

\begin{thm}\label{mainth} Let $n_1,...n_k$, $m_1,...,m_k$ be 
nonzero integers. Let $d=gcd(n_1,...,n_k)$.  
Then the identity 
\begin{equation}
\widehat{G_{\la_1,...\la_k}^{n_1,...,n_k,a}}
=CG_{\eta_1,...,\eta_k}^{m_1,...,m_k,b},\label{4.13}
\end{equation}
 between distributions on $F^k$ in the weak sense, is satisfied
if and only if 
\begin{equation}
\eta_i\la_i\nu_F=\gamma^{n_i}.\label{4.14}
\end{equation}
where $\gamma\in X(F^*)$ 
is a character, and one of the following two conditions holds: 

1. $\sum m_i=2$, $m_i=n_i$, and

\begin{equation}
\Gamma^F_{d,a}(u^d)\prod_{i=1}^k\Gamma^F(u^{-n_i}\la_i\nu_F)=
C\Gamma^F_{d,b}(u^{-d}\gamma^d)
;\label{4.15}
\end{equation}

2. $\sum m_i=0$, $m_i=-n_i$, and 
\begin{equation}
\Gamma_{d,a}^F(u^d)\prod_{i=1}^k\Gamma^F(u^{-n_i}\la_i\nu_F)=
C\Gamma_{d,b}^F(u^d\gamma^{-d}).
\label{4.17}
\end{equation}
\end{thm}

\Pf . Using the previous two lemmas, we get
\begin{equation}
\begin{align*}
&\int_{U(F^*)}\Gamma^F(u)u^{-1}(a)\prod_{i=1}^k\Gamma^F(\la_iu^{-n_i}\nu_F)
\la_i^{-1}u^{n_i}\nu_F^{-1}(p_i)du=\\
&C\int_{U(F^*)}\Gamma^F(v)v^{-1}(b)\eta_1v^{-m_1}(p_1)...\eta_kv^{-m_k}(p_k)dv.
\end{align*}\label{4.18}
\end{equation}
 
We see that this identity can hold only if the vectors
$(m_1,...,m_k)$ and $(n_1,...,n_k)$ are proportional. 
Let $\alpha=m_i/n_i$ for all $i$
(the proportionality coefficient). 
It is not difficult to see by asymptotic analysis of the above formula 
for $u,v=\nu_F^s$ for large $s$ (which is essentially 
equivalent to the ``semiclassical analysis'' of Section 2)
that $1-\sum m_i=\alpha$. Therefore, $1-\sum n_i=\alpha^{-1}$, so both 
$\alpha$ and $\alpha^{-1}$ are integers and hence $\alpha=\pm 1$. 
Thus, $m_i=\pm n_i$.

So we should consider two 
cases. 

{\bf Case 1.} $\sum m_i=2$, $m_i=n_i$. In this case,
replacing $v$ with $v^{-1}$, we obtain
\begin{equation}
\begin{align*}
&\int_{U(F^*)}\Gamma^F(u)u^{-1}(a)\prod_{i=1}^k\Gamma^F(\la_iu^{-n_i}\nu_F)
\la_i^{-1}u^{n_i}\nu_F^{-1}(p_i)du=\\
&C\int_{U(F^*)}\Gamma^F(v^{-1})v(b)\eta_1v^{n_1}(p_1)...\eta_kv^{n_k}(p_k)dv.
\end{align*}\label{4.18a}
\end{equation}

Let us replace $u$ with $u^{1/d}$. This leads to summation over
all roots of degree $d$, and therefore 
the Gamma functions $\Gamma^F(u)$, $\Gamma^F(v)$ are replaced
by the generalized Gamma functions:
\begin{equation}
\begin{align*}
&\int_{U(F^*)}\Gamma_{d,a}^F(u)\prod_{i=1}^k\Gamma^F(\la_iu^{-n_i/d}\nu_F)
\la_i^{-1}u^{n_i/d}\nu_F^{-1}(p_i)du=\\
&C\int_{U(F^*)}
\Gamma^F_{d,b}(v^{-1})\eta_1v^{n_1/d}(p_1)...\eta_kv^{n_k/d}(p_k)
dv.
\end{align*}\label{4.18b}
\end{equation}

It is clear that the integrals on the two sides of this equation
can coincide if and only if the contour of integration in the
first integral can be shifted to obtain the second integral. 
In particular, there must exist a character
$\gamma$ such that $\eta_i\la_i\nu_F=\gamma^{n_i}$. 
In this case, replacing $v$ with $u\gamma^{-d}$
(i.e., shifting the contour of integration), we get 
\begin{equation}
\begin{align*}
&\int_{U(F^*)}\Gamma_{d,a}^F(u)\prod_{i=1}^k\Gamma^F(\la_iu^{-n_i/d}\nu_F)
\la_i^{-1}u^{n_i/d}\nu_F^{-1}(p_i)du=\\
&C\int_{U(F^*)}\Gamma^F_{d,b}(u^{-1}\gamma^d)\la_1u^{n_1/d}\nu_F^{-1}(p_1)...
\la_ku^{n_k/d}\nu_F^{-1}(p_k)
du.
\end{align*}\label{4.18c}
\end{equation}

\begin{rem} We can shift the contour of integration without
worrying about residues, since our integral identities are
understood in the weak sense, while residual contributions 
are distributions supported on the coordinate hyperplanes, 
which by definition do not affect identities in the weak sense.   
\end{rem}

The latter condition is equivalent to 
\begin{equation}
\Gamma_{d,a}^F(u)\Gamma^F(\la_1u^{-n_1/d}\nu_F)...
\Gamma^F(\la_ku^{-n_k/d}\nu_F)=
C\Gamma^F_{d,b}(u^{-1}\gamma^d).
\end{equation}
Since $\Gamma_{d,a}$, by definition, vanishes away from the $d$-th powers, 
we can replace in this condition the variable $u$ by $u^d$ 
without changing the condition. This yields the equation
(\ref{4.17}) in the
theorem.

{\bf Case 2.} $\sum m_i=0$, $m_i=-n_i$. 
In this case, we see similarly to case 1 
that (\ref{4.13}) is equivalent to the combination 
of the two conditions from part 2 of the theorem.  

The theorem is proved. 
\qed

\subsection{ Gamma functions of prehomogeneous vector spaces.}

The construction of the previous section can in fact be further generalized  
to a more general setting of prehomogeneous vector spaces. 

For simplicity let $F$ have characteristic zero. 
Let $\underline{V}$ be a prehomogeneous vector space 
of dimension $M$ over $F$ for a reductive
group $\underline{G}$. We assume that $V=\underline{V}(F)$ has a dense
orbit under the action of the group of points $G=\underline{G}(F)$,
and that the same is true for
the dual prehomogeneous vector space $V^*$. 

Let $f$ be a relative invariant of $V$ generating its invariant
ring. Suppose it has degree $D$. Let $V^*$ be the dual space 
and $f_*$ the relative invariant of $V^*$, which is the 
multiplicative Legendre transform of $f$. 
Assume for simplicity 
that the stabilizer $G_0$ of $f$ in $G$ acts with finitely many orbits 
on the varieties $f=0$, $f_*=0$ in $V,V^*$. 

For any multiplicative character $\la$ of $F^*$, consider the
function $\la(f(x))$ on $V$. Since $F$ has characteristic zero, 
it is known that this function defines a distribution 
on $V$ which meromorphically depends on $\lambda$.   

It is easy to see that 
for generic $\lambda$ the Fourier transform 
$\widehat{\la\nu_F^{-\frac{M}{D}}(f)}(x)$
is proportional to $\la^{-1}(f_*(x))$.
(For generic $\lambda$, this is the unique, up 
to scaling, distribution of the correct homogeneity degree,
due to our assumption about the finiteness of the number of orbits).
This allows one to define the Gamma function $\Gamma^V(\la)$ of
$V$ to be the meromorphic function on $X(F^*)$ given by 
\begin{equation}
\widehat{\la\nu_F^{-\frac{M}{D}}(f)}(x)=\Gamma^V(\la)\la^{-1}(f_*(x)),
\label{4.4p}
\end{equation}
(whenever both sides
define distributions on $F$). 

 From the inversion formula for Fourier transform 
one gets the functional equation
\begin{equation}
\Gamma^V(\la)\Gamma^V(\nu_F^{M/D}\la^{-1})=\la(-1)^D.\label{4.5a}
\end{equation}

In fact, it was shown by Sato and Shintani \cite{SS} 
that if $F$ is archimedean, then these results are 
valid without the assumption of the finiteness of the number
of orbits. 

\subsection{ Identities with monomials on prehomogeneous vector 
spaces.}

Let $V_1,...,V_k$ be prehomogeneous vector spaces over $F$ 
as in the previous section, of dimensions $M_i$, and let 
$f_i$ be their relative
invariants, 
of degrees $D_i$. 
Let $V=\oplus V_i$. 

By a monomial on $V$ we will mean a rational function 
on $V$ of the form \linebreak
$f_1(x_1)^{n_1}...f_k(x_k)^{n_k}$, where $n_i$ are integers, 
and $x_i\in V_i$. 
Define the distribution on $V$ by the formula

\begin{equation}
{\mathcal G}^{n_1,..,n_k,a}_{\la_1,...,\la_k}(x_1,...,x_k)=
\psi(a\prod_{i=1}^kf_i(x_i)^{n_i})\la_1(f_1(x_1))...\la_k(f_k(x_k)).
\label{4.102}
\end{equation}
Here $a\in F^*$ is a parameter, and $\la_i\in X(F^*)$. 
Similar to the case of the usual monomials, this 
is a distribution for large real parts of $\lambda_i$
which meromorphically extends to generic $\la_i$, 
and a distribution in the weak sense for all $\la_i$. 

For this class of distributions, we have the following generalization 
of Theorem \ref{mainth}. 

\begin{thm}\label{prehomth} 
Let $n_1,...n_k$, $m_1,...,m_k$ be 
nonzero integers. Let $d=gcd(n_1,...,n_k)$.  
Then the identity 
\begin{equation}
\widehat{{\mathcal G}_{\la_1,...\la_k}^{n_1,...,n_k,a}}
=C{\mathcal G}_{\eta_1,...,\eta_k}^{m_1,...,m_k,b},\label{4.13b}
\end{equation}
 between distributions $V^*=\oplus V_i^*$ in the weak sense is satisfied
if 
and only if 
\begin{equation}
\eta_i\la_i\nu_F^{M_i/D_i}=\gamma^{n_i},\label{4.14b}
\end{equation}
where $\gamma\in X(F^*)$ 
is a character, and 
one of the following two conditions holds: 

1. $\sum m_iD_i=2$, $m_i=n_i$, 
\begin{equation}
\Gamma^F_{d,a}(u^d)\prod_{i=1}^k\Gamma^{V_i}
(u^{-n_i}\la_i\nu_F^{M_i/D_i})=
C\Gamma^F_{d,b}(u^{-d}\gamma^d)
;\label{4.15b}
\end{equation}

2. $\sum m_iD_i=0$, $m_i=-n_i$, and
\begin{equation}
\Gamma_{d,a}^F(u^d)\prod_{i=1}^k\Gamma^{V_i}(u^{-n_i}\la_i
\nu_F^{M_i/D_i})=
C\Gamma_{d,b}^F(u^d\gamma^{-d}).
\label{4.17b}
\end{equation}
\end{thm}

\begin{rem} Here the distribution 
${\mathcal G}_{\eta_1,...,\eta_k}^{m_1,...,m_k,b}$
is defined using the relative invariants $f_i^*$ of $V_i^*$, 
which are multiplicative Legendre transforms of $f_i$.
\end{rem}

\Pf . The proof of this theorem is analogous to the proof for
usual monomials. 
\qed

An important special case of this theorem describes integral identities 
with monomials of norms of field extensions, which includes the
identity from \cite{K1} cited in Section 1.1. This special 
case is defined by setting $V_i=F_i$ (field extensions of
$F$ of degrees $d_i$), and $f_i=\Nm:F_i\to F$ to be the norm maps.
In this case, 
$V_i$ is identified with $V_i^*$ using the trace functional.  

\subsection{ Poles of distributions
  $G_{\la_1,...,\la_k}^{n_1,...,n_k,a}$ }

We would like to find out when integral identities 
with monomials hold not only in the weak but also in 
the strong sense (i.e., on all test functions). 
For this purpose we should first of all find out 
where (in terms of $\la_i$) the distribution 
  $G_{\la_1,...,\la_k}^{n_1,...,n_k,a}$ 
could have poles. 
This is our goal in this section. 
For brevity we will denote this distribution just by $G$. 

\begin{prop}\label{poles} 
The divisor of poles of $G$ is a subset of the set of points 
where one of the following equations is satisfied:

(i) $\la_i$ is a singular character for some $i$ such that 
$n_i\ge 0$;

(ii) $(\la_jr_j^{-1})^{n_l}=(\la_lr_l^{-1})^{n_j}$
for some $(j,l)$ such that $n_j>0>n_l$,
and singular characters $r_j,r_l$.
\end{prop}

(The notion of a singular character was introduced in Section
4.2).

\Pf . It is clear that we can assume that 
all the exponents $n_i$ are nonzero. 

Suppose that a point $(\la_1,...,\la_k)\in X(F^*)^k$ is a 
generic point of the divisor of poles of $G$. Then the leading 
coefficient $G_{top}$ of the Laurent expansion of $G$ at that point
is a distribution on the coordinate cross. Its support is a closed subset
$S$ of the coordinate cross which is invariant under the scaling 
group $T=\{(t_1,...,t_k)\in (F^*)^k:\prod t_i^{n_i}=1\}$. 

Let $H_I=\{x_i=0,i\notin I\}$, $I\subset \{1,...,k\}$. 
Assume that $\la_i$ is not a singular character whenever 
$n_i>0$. Then  
$S\subset \cup_{I:|I|=k-2} H_I$.

Let $y$ be a generic point of $S$. In this case the $T$-orbit of
$y$ spans a certain subspace $H_I$, $|I|\le k-2$. In a small neighborhood of
$y$, the sets $S$ and $H_I$ coincide. 
 
Let $T_I$ be the group of elements of $T$ in which 
$t_i=1$, $i\in I$ (so $T_I$ acts trivially on $H_I$). 
On the one hand, for $t\in T_I$ we have 
$G_{top}(tx)=\prod_{i\notin I} \lambda_i(t_i)G_{top}(x)$. 
On the other hand, since 
in the neighborhood of $y$ the distribution $G_{top}$ is supported 
on $H_I$,
there exist singular characters $r_j$, $j\notin I$ 
such that $G_{top}(tx)=\prod_{i\notin I}r_i(t_i)G_{top}(x)$. 
This implies that whenever 
$\prod_{i\notin I}t_i^{n_i}=1$, we have 
$\prod_{i\notin I}(\la_ir_i^{-1})(t_i)=1$, for suitable 
singular characters $r_i$. In particular, for any distinct 
$i,j\notin I$, one has $(\la_ir_i^{-1})^{n_j}=
(\la_jr_j^{-1})^{n_i}$. 
 
The last equation can define a component of the pole divisor only if 
$n_j,n_l$ have opposite signs, since otherwise there are solutions 
with $\Re(\la_j)>0$, $\Re(\la_l)>0$, where $G$ is holomorphic. 
The proposition is proved. 
\qed

\subsection{ Identities with monomials in the weak and strong
  sense }

Sometimes one can deduce from an integral identity in the weak
sense that it actually holds in the strong sense (i.e., 
on all test functions). Let us do it in the case of archimedean
fields, using the theory of D-modules.

Let $F$ be archimedean. Let us say that a 
nonsingular multiplicative
character $\la$
is {\it strongly regular} if it generates 
(as a distribution) an irreducible D-module on the
line. For $F=\R$  a character is strongly regular iff 
it is not of the form $x^n/|x|$ for an integer $n$, 
or $x^n$ for $n<0$. 
For $F=\C$, a character is strongly regular iff 
it is not of the form $x^n\bar x^m$, where $n,m$ are
integers, and at least one of them is negative.  

\begin{thm} \label{weakstrong}
Suppose that for some k-tuples of characters 
$(\la_1,...,\la_k)$, $(\eta_1,...,\eta_k)$ 
identity (\ref{4.13}) holds in the weak sense, and 

(i) $\lambda_i$, $\eta_i$ are strongly regular whenever $n_i\ge 0$;

(ii) $(\la_jr_j^{-1})^{n_l}\ne (\la_lr_l^{-1})^{n_j}$
for any $(j,l)$ such that $n_j>0>n_l$,
and any singular characters $r_j,r_l$, and 
the same is true about $\eta_i$. 

In this case both sides of (\ref{4.13}) are well defined as
distributions (on all test functions), 
and (\ref{4.13}) is satisfied as an identity 
between distributions (in the ``strong'' sense). 
\end{thm}

\Pf . According to Proposition \ref{poles},
the distributions $G_1=G_{\la_1,...,\la_k}^{n_1,...,n_k,a}$, and 
$G_2=G_{\eta_1,...,\eta_k}^{m_1,...,m_k,b}$ are well defined. 
The fact that (\ref{4.13}) holds in the weak
sense means that $\hat G_1-G_2=\Delta$, where 
$\Delta$ is the sum of a distribution supported
on the coordinate cross and a distribution whose 
Fourier transform is supported on the cross.
Let $M_i$ be the D-modules generated by $G_i$.

{\bf Lemma.} $M_1$ and $M_2$ are irreducible. 

{\it Proof.}
Let us prove that $M_1$ is irreducible. 
The irreducibility of $M_2$ is shown in a similar way. 

Let $I\subset \{1,...,k\}$. We claim that if the numbers
$n_i,i\in I$, have the same signs, then the D-module $M_1$ is 
irreducible on the formal polydisk around 
a generic point $y$ of the subspace $H_I$ 
(where $H_I$ was defined in the proof of Proposition \ref{poles}).
Indeed, if all $n_i,i\in I$ are positive, this restriction 
is isomorphic to the exterior tensor product 
of $<\la_i>,i\in I$ with $k-|I|$ copies of 
$<1>$ (the structure sheaf of the 1-dimensional formal disk), 
so it is irreducible by assumption (i). 
If all $n_i$ are negative, the restriction 
is generated by $\prod_{j\in I}\la_j(x_j)
e^{i\text{Re}(c(x)\prod_{j\in I}x_j^{n_j})}$, 
where $c$ depends only of $x_j,j\notin I$, and is nonzero at $y$. 
This D-module is obviously irreducible for any $\la_j$. 

Thus, we see that all Jordan-Holder components of 
$M_1$ whose support is contained in the cross, 
are supported on the union 
of subspaces $H_{I_{j,l}}$ for $n_j>0>n_l$, where $I_{j,l}$ 
is the complement of $\{j,l\}$. 
Now, using scaling arguments as in the proof 
of Proposition \ref{poles}, it is easy to deduce from condition 
(ii) for $j,l$ that all these components are zero.
This means that $M_1$ is irreducible.  
The lemma is proved. \qed

Now let us prove the theorem. By the irreducibility of $M_1$,
 $\Delta$ is supported on the cross 
(as well as its Fourier transform), since the restriction of
$M_1$ to the set of its smooth points has to be an irreducible
local system. 

{\bf Lemma.}
The distribution $\Delta$ has support of codimension 2 or more.

{\it Proof.} If the support of $\Delta$ contains a point 
$(x_1,...,x_k)$ with $x_i=0$, $x_j\ne 0$ for $j\ne i$ then 
let us apply to $\hat G_1$ the algebra of differential operators 
in one variable $x_i$, with coefficients 
in polynomials of other variables. It suffices to check that 
the space of obtained distributions contains a 
nonzero distribution supported on the cross
(this would contradict the irreducibility of $M_2$). 

For this, it suffices to check that 
for any distribution on the line of the form
$h=\chi(t)\psi(ct^m)+\Delta$, 
where $m\ne 0$, $\chi$ is strongly regular if $m>0$, 
and $\Delta$ is supported at zero, $\Delta\ne 0$, 
one can find a polynomial differential operator
$D$ (algebraically depending on $c$) 
such that $Dh$ is nonzero but is supported at zero. 
This is straightforward. 
\qed 

Finally, since $\text{Supp}(\Delta)$ has codimension 2 or more, 
it is not difficult to see
using homogeneity arguments that $\Delta=0$. 

The theorem is proved.\qed

\begin{rem} The same method can be used for identities with norms
and prehomogeneous vector spaces; however, we do not give here
the details of such arguments. 
\end{rem}

\section{Identities over $\R$ and $\C$.}

As we saw in the previous section, 
the Gamma functions for $\R$ and $\C$ are given by simple
explicit formulas via the Euler Gamma function. One can also compute 
the Gamma functions of prehomogeneous 
vector spaces over these fields, using the theory 
 of Bernstein polynomials.
This allows one to express the Gamma function identities 
of the previous section in much more explicit (combinatorial)
terms, thus giving a completely elementary criterion of the
existence of integral identities. This is what we do in this section. 

For the sake of brevity, we do a systematic analysis only 
in the case of ordinary monomials with relatively prime
exponents, restricting ourselves to a number of examples in other
cases. 

\subsection{ The relation of the parameters $a$ and $b$.}

\begin{lem}
Suppose that identity (\ref{4.15b}) is satisfied in the weak
sense. Then if $\sum m_i=2$, one has 
\begin{equation}
ab=-\prod_i n_i^{-n_iD_i},\label{4.26}
\end{equation}
and if $\sum m_i=0$, one has
\begin{equation}
ab^{-1}=\prod_i n_i^{-n_iD_i}\label{4.26a}
\end{equation}
\end{lem}

\Pf . This follows from Theorem \ref{classical} or from analyzing
the asymptotics of the Gamma function identity for 
large $s$. 
\qed

\subsection{  Multiplication formulas for $\Gamma$-functions.}
 
We recall the classical multiplication law for $\Gamma$-function:
\begin{equation}
\Gamma(Nz)=
N^{Nz-\frac{1}{2}}
(2\pi)^{(1-N)/2}\Gamma(z)\Gamma(z+1/N)...\Gamma(z+(N-1)/N).\label{4.21}
\end{equation}

This implies that
\begin{equation}
\Gamma^\C_{Nn}(Ns)=N^{Ns-1}\Gamma^\C_n(s)\Gamma^\C_n(s+2/N)...
\Gamma^\C_n(s+2(N-1)/N),\label{4.22}
\end{equation}
and if $N$ is odd then 
\begin{equation}
\Gamma^\R_{Nn}(Ns)=i^{-(N-2)(N-1)/2}
N^{Ns-1/2}\Gamma^\R_n(s)\Gamma^\R_{n+1}(s+1/N)...
\Gamma^\R_{n+N-1}(s+(N-1)/N),\label{4.23}
\end{equation}

This can be uniformly written as 
\begin{equation}
\Gamma^F(\chi^N)=C_{F,N}
\chi(N^N)\Gamma^F(\chi)\Gamma^F(\chi\text{Nm}^{1/N})...
\Gamma^F(\chi\text{Nm}^{(N-1)/N}),\label{4.24}
\end{equation}
($N$ is odd for $F=\R$), where $\text{Nm}$ is the algebraic norm of $F$ 
as an extension of $\R$ ($\Nm(x)=x\bar x$ for $\C$ and $\Nm(x)=x$
for $\R$), and $C_{F,N}=N^{-1}$ for $\C$ and
$i^{-(N-2)(N-1)/2}N^{-1/2}$
for $\R$.   

\subsection{  Classification of identities with monomials.}

Define the group 
$X_F=X(F^*)/<\text{Nm}>$. Let $Div(X_F)$ be the group of divisors
on $X_F$. Now for $\chi\in X_F$ and a positive integer $N$ set 
$D_{\chi,N}=(\chi)+(\chi\text{Nm}^{1/N})+...+(\chi\text{Nm}^{(N-1)/N})
\in Div(X_F)$
(here $N$ should be odd if $F=\R$). 
Set $D_{\chi,-N}=-D_{\chi^{-1}\nu_F,N}$. 

\begin{thm} \label{divrel}
(i) Make the assumptions of Theorem \ref{mainth},
and suppose that the exponents $n_i$ are relatively prime: $d=1$.  
Then identity (\ref{4.13}) holds in the weak sense for some $C$, 
$a,b,\lambda_i,\eta_i$
in and only in one of the following two situations:

1. $\sum m_i=2$, $m_i=n_i$, 
and there exists $\xi,\mu_i\in X_F$ such that:
\begin{equation}
D_{1,1}+\sum_i D_{\mu_i,-n_i}=D_{\xi,-1}.\label{4.25}
\end{equation}

2. $\sum m_i=0$, $m_i=-n_i$, 
and there exists $\xi,\mu_i\in X_F$ such that:
\begin{equation}
D_{1,1}+\sum_i D_{\mu_i,-n_i}=D_{\xi,1}.\label{4.27}
\end{equation} 

(ii) More precisely, 
in any of the above two cases, one can find $\lambda_i,\gamma\in
X(F^*)$, $a,b\in F^*$ 
(where $\gamma$ is as in Theorem \ref{mainth}), 
so that the integral identity of Theorem \ref{mainth} is
satisfied, and  
$\lambda_i\nu_F=\mu_i$, $\gamma=\xi$ (resp. $\gamma=\xi^{-1}$) modulo $\Nm$
in the first (resp. second) case.
\end{thm}

\Pf . Let us sketch the proof of (i). The proof of (ii) is
obtained automatically in the process of proving (i). 

To prove the necessity of the conditions of (i), let us
represent the elements of the group $X_F$ in $X(F^*)$ 
by a fundamental domain lying in the region $\text{Re} u>>0$. 
Let us write down the identity of divisors of zeros and poles 
for the left and right
hand sides of the two Gamma function relations of Theorem \ref{mainth}
in this fundamental domain. It is easy to check that this yields 
exactly the two divisorial relations above. 

Let us now prove the sufficiency of (i). It is not difficult to show that if 
the conditions of (i) are satisfied, then one can find elements
$\lambda_i,\gamma\in X(F^*)$ in the cosets 
$\mu_i\nu_F^{-1},\xi^{\pm 1}$ of the free cyclic group $<\Nm>$, such that 
the divisors of both sides of the Gamma function relations of
Theorem \ref{mainth} in $X(F^*)$ (not only in $X_F$) coincide. 

By the multiplication formula (\ref{4.24}),
this implies that these relations hold up to a constant and an 
exponential factor. These factors can be removed by 
an appropriate choice of the constant $C$, and by 
imposing relations (\ref{4.26}),(\ref{4.26a}) on $a$ and $b$. 
\qed

\subsection{ Identities with monomials and 
exact covering systems.}

Let $(p_1,...,p_n)$ be positive integers whose sum equals $p$. 
According to a classical definition in elementary number theory
(\cite{G}, problem F14) 
an exact covering system 
of type $(p_1,...,p_n)$ is a representation of the group 
$\Z/p\Z$ as a disjoint union of cosets of 
 groups $\Z/p_i\Z$, 
$i=1,...,n$. 
It is clear that for the existence of an exact 
covering system of type $(p_1,...,p_n)$, it is
necessary that all $p_i$ are divisors of $p$. But this is not
sufficient: for example, it is clear that for a covering
system, $p/p_i$ is never relatively prime to $p/p_j$. 
In fact, there is no known simple necessary and sufficient condition 
for the existence of an exact covering system,
and various questions about exact covering systems
are the subject of an extensive theory 
(see \cite{G}, problem F14, and references therein). 
It is interesting, therefore, that in a special case, 
integral identities with monomials that we considered 
correspond to exact covering systems. 

Consider integral identities of type (\ref{4.13})
with $\sum n_i=2$, $n_k=n>0$, and $n_i<0$ for $i=1,...,k-1$. 
If $F=\R$, we assume that $n$ is odd (this is clearly necessary 
for the existence of integral identities, since the multiplication
formula over $\R$ is valid only for multiplication by an odd
number). 

We will realize
$\Z/n\Z$ as the set of integers $\{0,...,n-1\}$ and parametrize
covering systems by sequences $(p_1,...,p_{k+1})$ where 
$p_i$ is the biggest element of the i-th coset of the covering. 

\begin{thm}\label{covsys} Integral identities of type
  (\ref{4.13})
 with such data
for $F=\C$ are (up to rescaling $a$ and $b$) 
in 1-1 correspondence with exact covering systems 
of type $(-n_1,...,-n_{k-1},1,1)$. Moreover, 
the formula corresponding to the covering system 
$(p_1,...,p_{k+1})$ has the following parameters:
\begin{equation}
\la_k=\Nm^{n-1-p_k},\
\la_j=\nu_F^{-1}\Nm^{\frac{(p_j-p_k)n_j}{n}},
\end{equation}
\begin{equation}
\eta_k=\Nm^{p_{k+1}},
\eta_j=\nu_F^{-1}\Nm^{\frac{(p_{k+1}-p_j)n_j}{n}}.
\end{equation}
\end{thm}

\Pf . It is clear from the divisorial relation of 
Theorem \ref{divrel}
that to any integral identity with the given data, there
corresponds a canonical exact covering system. Indeed, consider
the divisor $D_{\mu_k^{-1}\nu_F,n}$. According to Theorem
\ref{divrel}, it is exactly covered by the divisors
$D_{\mu_i,-n_i}$, $i=1,..,k-1$, $D_{1,1}$, and 
$D_{\chi^{-1},1}$. Let us regard
the divisor $D_{\mu_k^{-1}\nu_F,n}$ as
$\Z/n\Z$ by declaring $\mu_k^{-1}\nu_F$ to be the unit. 
Then the other divisors define an exact covering of the group. 

Conversely, given an exact covering system, the characters
$\la_i$, $i<k$, are uniquely determined by $\la_k$
from the condition of cancellation of zeros and poles. 
Moreover, $\la_k$ itself is uniquely determined, 
because the Gamma-factor corresponding to the divisor $D_{1,1}$ 
does not have a shift by a character. 
In fact, a direct computation shows that 
$\la_i$ are given by the explicit formulas in the theorem.
The same computation shows that
$\gamma=\nu_F\Nm^{(p_{k+1}-p_{k})/n}$, 
which allows one to compute $\eta_i$ from the formula 
$\la_i\eta_i\nu_F=\gamma^{n_i}$.  
The theorem is proved. 
\qed

\begin{thm} Integral identities corresponding to 
exact covering systems hold not only in the weak sense, 
but also in the strong sense. 
\end{thm}

\Pf . It is clear from Theorem \ref{covsys} that 
$\la_k$ is a monomial, so it is always strongly regular. 
Therefore, by Theorem \ref{weakstrong}, it suffices to check 
that there is no $j\le k-1$ and singular characters $r_j,r_k$ 
such that $(\la_jr_j^{-1})^n=(\la_kr_k^{-1})^{n_j}$. 
Suppose such $j,r_j,r_k$ exist. Using the explicit formulas
for $\la_i$, it is easy to deduce from this that 
$2+p_j\ge n(1-1/n_j)$, which is impossible. 
The theorem is proved.
\qed

\subsection{An  example}

Consider the simplest example of the above theory: 
$k=2$, $n=n_2=3$, $n_1=-1$, $F=\R$.  
Applying the theorems we have proved, we obtain the following
result. 

\begin{thm}
Consider the distribution
\begin{equation}
\psi(ix^3/y)\la_2(x)\la_1(y).\label{4.28}
\end{equation}

(i) Such a distribution has an elementary Fourier transform 
in the weak sense if and only if it is one of the following 
six distributions:
\begin{equation}
\begin{align*}
&G_1=e^{ix^3/y}|y|^{-1/3};
G_2=e^{ix^3/y}|y|^{-2/3}\text{sign}(y);
G_3=e^{ix^3/y}x|y|^{-4/3}\text{sign}(y);\\
&G_4=e^{ix^3/y}x^2|y|^{-5/3};
G_5=e^{ix^3/y}x|y|^{-2/3}\text{sign}(y);
G_6=e^{ix^3/y}x^2|y|^{-4/3}\text{sign}(y);
\end{align*}
\end{equation}

(ii) The Fourier transform acts on these distributions by 
$$
\hat G_j(x,y)=\pm i G_{s(j)}(x,-y/27), 
$$
where the sign is $+$ for $j=2,3,5,6$, and $-$ for $1,4$. and $s$  
is the involution of $\{1,2,3,4,5,6\}$ given by $s=(13)(45)$.
These identities hold in the strong sense.  
\end{thm}

\subsection{ Calculation of Gamma functions 
of prehomogeneous vector spaces over $\R$ and $\C$.} 

Now we want to study integral identities for prehomogeneous
vector spaces. For this purpose, we need an explicit 
expression for Gamma functions of these spaces. 
Luckily, 
if $F=\R$ or $\C$ then the
Gamma function (up to a constant) can be found in a number of
cases using
Bernstein's polynomial. This is well known
(see e.g. \cite{SS}), but we will give the
argument for the reader's convenience.  

For simplicity we restrict ourselves to the case $F=\C$. 

Consider the (multivalued) 
function $f_*(\dd)f^{s+1}$ for $s\in\C$. It is easy to see
that there exists a unique monic polynomial ${\bold b}(s)$ of degree $D$ 
such that $f_*(\dd)f^{s+1}={\bold b}(s)f^{s}$
(Bernstein's polynomial). In the case of prehomogeneous vector
spaces that we are considering, 
Bernstein's polynomial was introduced by Sato in the sixties and 
is called ``Sato's b-function''
(\cite{Sa}). For irreducible regular spaces, the b-functions
were computed by Kimura in \cite{Ki}. 

One has (for generic $s$):
\begin{equation}
f_*(\dd)\la_{s+2,n}(f)={\bold b}(\frac{s-n}{2})\la_{s+1,n-1}(f),
\overline{f_*(\dd)}\la_{s+2,n}(f)={\bold b}(\frac{s+n}{2})\la_{s+1,n+1}(f).
\label{diffeq}
\end{equation}

Now let us find the Gamma function.
We have 
$$
\widehat{\overline{f_*(\dd)}h}(x)=2^{-D}(-i)^D{f_*(x)}\hat h(x),
\widehat{{f_*(\dd)}h}(x)=2^{-D}(-i)^D\overline{f_*(x)}\hat h(x).
$$

Therefore, we have
$$
{\bold b}(\frac{s\mp n}{2}-\frac{M}{D})
\Gamma^V(\la_{s-1,n\pm 1})=2^{-D}(-i)^D\Gamma^V(\la_{s,n}).
$$
 
\begin{rem} The last equation and the functional equations
for the Gamma functions
imply that 
$$
{\bold b}(s-M/D)=(-1)^D{\bold b}(1-s).
$$
This shows that the collection of roots of the polynomial ${\bold b}(s)$ is symmetric
with respect to the point $-\frac{1}{2}(1+\frac{M}{D})$. 
\end{rem}

 From now on we will assume that 
$M/D$ is an integer, and that the roots of ${\bold b}(s)$ are
 integers or half integers. This condition is satisfied for a number 
of cases in \cite{Ki}.  
In this case, the obtained difference equation 
 allows one to deduce a formula for the Gamma function
(up to a scalar). Namely, we have

\begin{prop} \label{gammapreh}
$$
\Gamma^V(\la_{s,n})=C_V\prod_{j=1}^D\Gamma^F(\la_{s+2(s_j-M/D),n}),
\label{gammav}
$$
where ${\bold b}(s)=\prod_i(s+s_i)$, and $C_V$ is a constant. 
\end{prop}

\Pf . 
Denote the proportionality coefficient 
between the LHS and the RHS of (\ref{gammav}) by 
$C_V(s,n)$. It is obvious that this function is periodic: 
$C_V(s+1,n\pm 1)=C_V(s,n)$ and satisfies
$C_V(n,s)C_V(-n,\frac{2M}{D}-s)=1$
(by virtue of the functional equation).
Thus, it suffices to show that $C_V(s,0)$ is constant. 
 
Since $\Gamma^V(\la_{s-2,0})=c\cdot
{\bold b}(1-s)^{-2}\Gamma^V(\la_{s,0})$,
and $\Gamma^V(\la_{s,0})$ is holomorphic for large 
positive $\text{Re}(s)$, 
the poles of $\Gamma^V(\la_{s,0})$ 
can only arise at points of the form 
$2(s_j+r)$, where $r$ is an integer.
Since $s_j$ are integers
or half integers, this means that all the poles are
integers. 
Similarly, because $1/\Gamma^V$ is holomorphic for large 
negative $\text{Re}(s)$,  
we find that all the zeros are integers. 
The same clearly holds for the product of the Gamma functions 
on the right hand side of (\ref{gammav}). 
Therefore, the same holds for $C_V(s,0)$. 

Thus, we have: the meromorphic function $h(s):=C_V(s,0)$ 
is periodic with period 2, its zeros and poles are integers, and 
$C_V(s,0)C_V(-s,0)=1$ (since $2M/D$ is an even integer). 
This implies that the function must have no zeros and no poles, since
any zero of this function will also be its pole, and vice versa. 

Given this, an asymtotic analysis for large $s$ 
(using the stationary phase approximation)
shows that $C_V(n,s)$ is a constant. 
Another way to show it (knowing that 
$C_V$ has no zeros or poles) is to use the result of 
\cite{SS}, which shows that 
$C_V$ is a trigonometric function.  
\qed

\subsection{ Identities for prehomogeneous vector spaces }

Let $W$ be a prehomogeneous vector space of dimension $M$ 
over $F=\C$, satisfying the
conditions that we imposed in the previous section. Let $f$ be its 
relative invariant, of degree $D$. Consider the space
$V=V_1\oplus...
\oplus V_{D-1}$, where $V_{D-1}=W$ and $V_j=F$ for
$j=1,...,D-2$. We will denote elements of $W$ by $x$ and 
elements of $V_j$ by $t_j$ for $j=1,...,D-2$. 
Let $f_j=t_j$, $j=1,...,D-2$, and $f_{D-1}=f$. 
Let $n_1=...=n_{D-2}=-1$ and $n_{D-1}=D$. 
We are interested in identities of type (\ref{4.13b})
arising in this situation. 

\begin{thm} Identities of type (\ref{4.13b}) 
for the described data, up to rescaling $a$ and $b$, 
correspond to functions $\sigma:\{1,...,D\}\to \Bbb Q$ 
such that $\prod_j(s+\sigma(j))$ equals the b-function 
${\bold b}(s)$ of the space $W$. More precisely, for any such function
$\sigma$ the parameters $\la_j,\eta_j$ of 
the corresponding identity (\ref{4.13b}) are given by 
the formula
\begin{equation}
\la_i=\nu_F^{\sigma(i)-\sigma(D-1)-1}, 
i\le D-2,\
\la_{D-1}=\nu_F^{\sigma(D-1)-M/D},
\end{equation}
\begin{equation}
\eta_i=\nu_F^{\sigma(D)-\sigma(j)-1}, i\le D-2,\
\eta_{D-1}=\nu_F^{1-\sigma(D)}.
\end{equation}
\end{thm}

\Pf . The proof is analogous to the proof for usual monomials. 
\qed

\subsection {Example.} 
Let $W$ be the 27-dimensional prehomogeneous
vector space over $\C$ from example 6 in Section 3.4. 
Let $W^*$ be the dual space, and $f_*$ the 
multiplicative Legendre transform of $f$. 
The b-function of this space equals $(s+1)(s+5)(s+9)$ (see 
\cite{Ki}). Therefore, from the previous theorem we 
obtain the following result.

\begin{thm}
Consider the distribution $\psi(i\text{Re}
f(x)/y)|f(x)|^p|y|^q$ on $W\oplus \C$ ($p,q\in \C$). 

(i) Such a distribution has an elementary Fourier 
transform in the weak sense if and only if 
$(q,p)$ takes one of the following six values:
$(-10,-8); (-18,0); (6,-16);$ $(-10,0);
(6,-8); (14,-16)$.

(ii) The Fourier transform in the weak sense maps 
any of these distributions to a distribution of the same type
on $W^*\oplus \C$ (up to a scalar). The action of 
the Fourier transform on the parameters $(q,p)$ is given by the
 following involution:
$$
(-10,-8)\to (14,-16),\ (-18,0)\to (6,-8),\
(6,-16)\to (6,-16),
$$
$$ (-10,0)\to (-10,0).
$$

The identity corresponding to 
$(-10,0)\to (-10,0)$ holds in the strong sense. 
\end{thm}

\begin{rem} The fact that the identity 
corresponding to 
$(-10,0)\to (-10,0)$ holds in the strong sense is proved
similarly to the case of usual monomials. 
\end{rem}

\begin{rem} Note that the integral
identities 
other than $(-10,0)\to (-10,0)$ cannot 
be understood in the strong sense,  
because the distributions we considered have poles 
on the divisor $f(x)=0$ and therefore are not well defined 
without a regularization. So to replace the identities 
we considered by identities in the strong sense one would 
first need to regularize one or both sides. 
We leave this beyond the scope of this paper.  
\end{rem}

\section{Identities over local non-archimedean fields}

\subsection{Local constants and Gamma functions}
Let $F$ be a local non-archimedean field.
Let $\O\subset F$ be the valuation ring, 
$\pi\in F^*$ be a uniformizing element.
A character $\chi\in X(F^*)$ is called unramified if $\chi|_{\O^*}=1$. 
The local $L$-factors are defined as the following meromorphic functions of 
$\chi\in X(F^*)$: 
$$L(\chi)=\cases 1 &\ \text{ if }\chi \text{ is ramified}, \\
(1-\chi(\pi))^{-1} &\ \text{ if }\chi \text{ is unramified}.\endcases
$$
Note that $L(\chi)$ has a unique pole at the trivial character.
Then we have
$$\Ga^F(\chi)=\frac{L(\chi)}{L(\chi^{-1}\nu_F)}
\eps(\chi^{-1}\nu_F)$$ 
where $\eps(\chi)=\eps(\chi,\psi)\in\C^*$ are the local constants 
considered by Deligne
in \cite{De2}. With our choice of the Haar measure $dx$ on $F$ we have
$\eps(\chi)=1$ if $\chi$ is unramified.
If $\chi$ is ramified then
$$\eps(\chi,\psi)=\int_{F^*}\chi^{-1}(x)\psi(x)dx:=
\sum_n \int_{v(x)=n}\chi^{-1}(x)\psi(x)dx.$$
Note that the functional equation (\ref{4.5}) for Gamma function implies that
$$\eps(\chi)\eps(\chi^{-1}\nu_F)=\chi(-1).$$

\subsection{Identities between local constants}
Recall that for a local field $F$ the Weil group $W_F$ is defined
as the preimage of the subgroup generated by the Frobenius under
the surjective homomorphism $\Gal(\ov{F}/F)\ra\Gal(\ov{k}/k)$
where $k$ is the residue field of $F$.
The local class field theory provides an isomorphism between
the abelianization of the Weil group $W_F$ and $F^*$.
Following \cite{De2} we normalize this isomorphism in such a way
that uniformizing elements in $F^*$ correspond to liftings of the inverse of
the Frobenius. Note that for every finite separable extension $F\subset E$
we have a commutative diagram
\begin{equation}
\begin{array}{ccc}
W_E &\lrar{}& E^*\\
\ldar{i} & &\ldar{\Nm_{E/F}}\\
W_F &\lrar{}& F^*
\end{array}
\end{equation}
where $i$ is the natural inclusion.
Thus, to a character $\chi$ of $F^*$ we can associate a one-dimensional
representation of $W_F$. We denote by $[\chi]$ the corresponding element
of the representation ring $R(W_F)$ (of finite-dimensional 
continuous complex representations of $W_F$). If
$\la$ is a character of $F^*$ and $F\subset E$ is a finite separable extension
then 
$[\la\circ\Nm_{E/F}]=\Res[\la]$ where $\Res:R(W_F)\ra R(W_E)$ is the restriction
homomorphism.

The principal theorem of \cite{De2} (Theorem 4.1) says that the map
$[\chi]\mapsto\eps(\chi)$ extends to a homomorphism 
$V\mapsto\eps(V)$ from $R(W_F)$ to $\C^*$, such that
for a finite separable extension $E$ of $F$ and for any virtual representation
$V$ of $W_E$ of dimension $0$ one has
$$\eps(\Ind_{W_E}^{W_F}V,\psi)=\eps(V,\psi\circ\Tr_{E/F})$$
(where both sides do not depend on a choice of Haar measures).
If the dimension of $V$ is not zero then we can write
$$\eps(\Ind_{W_E}^{W_F}V)=\la(E/F)^{\dim V}\cdot\eps(V)$$
where $\la(E/F)=\eps(\Ind_{W_E}^{W_F}[\nu_F^s])$ (which does 
not depend on $s$).

On the other hand, one can extend local $L$-factors $L(\chi)$
to a homomorphism from
$R(W_F)$ to the group of non-zero meromorphic functions on $X(F^*)$ by setting 
$$L(V,\la)=\det(1-\Frob^{-1}|_{(\la\otimes V)^I})^{-1}$$
where $I\subset W_F$ is the inertia subgroup (see \cite{De2}, 3.5).
Moreover, for a finite separable extension $F\subset E$ and $[V]\in R(W_E)$ one has
$$L(\Ind_{W_E}^{W_F}[V],\la)=L([V],\la\circ \Nm_{E/F})$$
(see \cite{De2}, Prop. 3.8).

\begin{prop}
Let $F_i$ be finite separable extensions of a local field $F$,
and for every $i$ let $\chi_i$ be a character of $F_i^*$.
Assume that we have a linear relation 
$$\sum_i n_i\Ind[\chi_i]=0$$
between the induced representations $\Ind[\chi_i]$
of $W_F$, where $n_i\in\Z$.
Then one has the following identity:
$$\prod_i (\la(F_i/F)\Ga^{F_i}(\chi_i(\la\circ\Nm_{F_i/F}))^{n_i})=1$$
for $\la\in X(F^*)$.
\end{prop}

\Pf . It suffices to notice that we also have the linear relation
$$\sum_i n_i\Ind[\chi_i^{-1}\nu_{F_i}]=0$$
(since $\nu_{F_i}=\nu_F\circ\Nm_{F_i/F}$)
and apply the above properties of $L$-factors and $\eps$-constants.
\ed

Here is an example of an identity provided by the above proposition.
Let $F\sub E$ be a finite cyclic extension. 
Then for any character $\la$ of $F^*$ we have the relation
$$\Ind_{W_E}^{W_F}[\la\circ\Nm_{E/F}]=
\sum_{\chi\circ\Nm_{E/F}=1} [\chi\cdot\la]$$
where the sum is taken over all characters $\chi$ of $F^*$ which
are trivial on $\Nm_{E/F}(E^*)$ (note that this subgroup has index
$[E:F]$ in $F^*$).
Hence we derive the identity
$$\la(E/F)\cdot\Ga^E(\la\circ\Nm_{E/F})=\prod_{\chi\circ\Nm_{E/F}=1}
\Ga^F(\chi\cdot\la).$$ 

\begin{rem} In the simplest case $[E:F]=2$ 
this formula appears in \cite{GGP}.
\end{rem}

More generally, let $F\sub E$ be a finite cyclic extension.
For every $n>0$ such that
$n|[E:F]$ let $F_n\subset E$ be the cyclic
subextension of degree $n$ over $F$. 
Let us denote by $X_{F,E}$ the subgroup in $X(F^*)$ consisting
of characters of the form $\nu_F^m\chi$ where $m\in\Z$,
$\chi\in X(F^*)$ is trivial on $\Nm_{E/F}(E^*)\subset F^*$.
Let $Div(X_{F,E})$ be the group of divisors on $X_{F,E}$.
For every character $\chi$ in $X_{F_n,E}$
let us define the divisor
$D(\chi)\in Div(X_{F,E})$ 
as follows:
$$D(\chi)=\sum_{\la:\la\circ\Nm_{F_n/F}=\chi} (\la).$$
Note that the homomorphism
$$X_{F,E}\ra X_{F_n,E}:\la\mapsto \la\circ\Nm_{F_n/F}$$
is surjective, hence the divisor $D(\chi)$ has degree $n$.
As an additive character on $F_n$ let us take $\psi_n=\psi\circ\Tr_{F_n/F}$.
Then extending the Gamma function
to divisors multiplicatively we can write
$$\la(F_n/F)\cdot\Ga^{F_n}(\chi(\la\circ\Nm_{F_n/F}))=
\Ga^F(\la D(\chi)).$$
for $\chi\in X_{F_n,E}$, $\la\in X(F^*)$, where for every divisor $D$ we
denote by $\la D$ the divisor obtained from $D$ by applying
the shift $\mu\mapsto \la\mu$.

Now let $(d_1,\ldots,d_k)$ be a sequence of positive numbers
dividing the degree $[E:F]$, and let 
$(m_1,\ldots,m_k)$ be a sequence of integers.
For every $i=1,\ldots,k$ let $\chi_i$ be a character in $X_{F_{d_i},E}$. 
Assume that
$$\sum_i m_i D(\chi_i)=0.$$  
Then for every $\la\in X(F^*)$ we have
\begin{equation}
\prod_i \la(F_{d_i}/F)^{m_i}
\Ga^F(\chi_i(\la\circ\Nm_{F_{d_i}/F}))^{m_i}=1.
\end{equation}

Using the functional equation for Gamma functions we can rewrite this result 
slightly differently.
For $e=\pm 1$ and a character $\chi\in X_{F_n,E}$ let us set
$$D(\chi,e)=\cases D(\chi), & \text{ if } e=1,\\ 
-D(\chi^{-1}\nu_{F_n}), & \text{ if } e=-1 .\endcases$$
Then for a sequence $(e_1,\ldots, e_k)$, where $e_i=\pm 1$,
the relation
$$\sum_i D(\chi_i,e_i)=0$$
implies the identity
$$\prod_i \Ga^F(\chi_i(\la^{e_i}\circ\Nm_{F_{d_i}/F}))=C\cdot\la(\prod_i e_i)$$
for $\la\in X(F^*)$,
where the constant $C\in\C^*$ doesn't depend on $\la$.

\subsection{Identities for cyclic extensions}
Let $F\sub E$ be a cyclic extension of local fields,
$F_i\sub E$ be a subextension of degree $d_i$ over $F$ ($i=1,\ldots,k$). For
every $i=1,\ldots,k$ let
$\chi_i$ be a character in $X_{F_i,E}$, and let $e_i$ be either $1$ or $-1$.
Below we denote by $1_F$ the trivial character of $F^*$.
The following theorem follows easily from the above identities
for Gamma functions and from Theorem \ref{prehomth}.

\begin{thm} In the above situation
assume that one of the following identities in $Div(X_{F,E})$ holds:

1. $D(1_F,1)+\sum_i D(\chi_i,-e_i)=D(\xi,-1)$,

2. $D(1_F,1)+\sum_i D(\chi_i,-e_i)=D(\xi,1)$,\\
for some $\xi\in X_{F,E}$.
Then we have the following identity between distributions on
$\oplus_i F_i$ in the weak sense:
$${\mathcal F}(\prod_i(\chi_i\nu_{F_i}^{-1})(x_i)
\psi(a\prod_i\Nm_{F_i/F}(x_i)^{e_i}))=
C\cdot\prod_i\eta_i(x_i)\psi(b\prod_i\Nm_{F_i/F}(x_i)^{ee_i})$$
for some constants $C\in\C^*$, $a,b\in F^*$, where
$e=1$ in the case 1, $e=-1$ in the case 2,
$$\eta_i=\chi_i^{-1}(\xi^{ee_i}\circ\Nm_{F_i/F}).$$
\end{thm}

For example, the identity (\ref{cub}) corresponds to the following equality
of divisors:
$$D(1_F,1)+D({\mathcal E},1)+D(\nu_E,-1)=D({\mathcal E}\nu_F,-1).$$

Thus, we see that identity (\ref{cub}) holds in the weak sense. 

Let us prove now that identity (\ref{cub}) in fact holds in the
strong sense. We will give a sketch of the argument. 
First of all, the function $\phi_{{\mathcal E}}$
makes sense as a distribution (i.e., does not need
regularization). Since {\ref{cub}} holds in the weak sense, 
we have $\widehat{\phi_{\mathcal E}}=\epsilon \phi_{\mathcal E}
+\eta$, where $\eta$ is a distribution which is a sum 
of a distribution concentrated on the coordinate cross 
and a distribution whose Fourier transform is concentrated 
on the cross, and $\epsilon$ is a sign. 
The distribution $\eta$ has the same homogeneity 
properties as $\phi_{\mathcal E}$, and satisfies 
$\hat \eta=-\epsilon \eta$. From this it is easy to deduce that 
$\eta=c(\delta(t)-\epsilon \delta(x))$. Thus, it remains to check 
that $c=0$, which can be checked by a direct calculation;
for instance, one can check directly that 
$\hat \phi_{\mathcal E}$ is locally constant at $x=0,t\ne 0$.
Thus, (\ref{cub}) holds in the strong sense. 

\subsection{One more identity}
There are more complicated examples of identities between local constants
which involve non-abelian extensions (see \cite{De2}, section 1). Here is an example.
Let $l$ be a positive integer,
$E$ a Galois extension of $F$ with Galois group $G$ which is
a central extension of $(\Z/l\Z)^2$ by an abelian group $Z$.
Let $H_1\sub G$ (resp. $H_2$) be the preimages of the first (resp. the second)
factor $\Z/l\Z$ in $(\Z/l\Z)^2$. Let $\chi_i$ be a character of $H_i$
for $i=1,2$, such that $\chi_1|_Z=\chi_2|_Z$ and this character is non-trivial
on $[G,G]\sub Z$. Let $F_i\sub E$ be the subextension of $F$ corresponding
to the subgroup $H_i$. Then for $i=1,2$ we can consider $\chi_i$ as a character
in $X_{F_i,E}$ and for every character $\la\in X(F^*)$ one has
$$\Ind_{W_{F_1}}^{W_F} \chi_1(\la\circ\Nm_{F_1/F})=
\Ind_{W_{F_2}}^{W_F} \chi_2(\la\circ\Nm_{F_2/F})$$ 
Hence,
\begin{equation}
\la(F_1/F)\Ga^{F_1}(\chi_1(\la\circ\Nm_{F_1/F}))=
\la(F_2/F)\Ga^{F_2}(\chi_2(\la\circ\Nm_{F_2/F}))
\end{equation}
for $\la\in X(F^*)$.
Now the following theorem follows easily from Theorem \ref{prehomth}.

\begin{thm} In the above situation we have the following identity
of distributions on $F_1\times F_2$ in the weak sense:
$${\mathcal F}((\chi_1\nu_{F_1}^{-1})(x_1)\chi_2^{-1}(x_2)
\psi(\frac{\Nm_{F_1/F}(x_1)}
{\Nm_{F_2/F}(x_2)}))=
C\cdot(\chi_2\nu_{F_2}^{-1})(x_2)
\chi_1^{-1}(1_2)\psi((-1)^l\frac{\Nm_{F_2/F}(x_2)}
{\Nm_{F_1/F}(x_1)}).$$
\end{thm}

\end{document}